\documentclass[11pt]{article}

\usepackage{latexsym}
\usepackage[dvips]{graphicx}
\usepackage{float}
\usepackage{amsmath}
\usepackage{amssymb}
\usepackage{epic}
\usepackage{color}
\usepackage{setspace}
\usepackage{lscape}
\usepackage{arydshln}
\usepackage{fancybox}   
\usepackage{mathtools}
\usepackage{amsthm}

\DeclarePairedDelimiter\ceil{\lceil}{\rceil}

\begin{document}
\bibliographystyle{plain}
\floatplacement{table}{H}
\newtheorem{definition}{Definition}[section]
\newtheorem{lemma}{Lemma}[section]
\newtheorem{theorem}{Theorem}[section]
\newtheorem{corollary}{Corollary}[section]
\newtheorem{proposition}{Proposition}[section]

\newcommand{\sni}{\sum_{i=1}^{n}}
\newcommand{\snj}{\sum_{j=1}^{n}}
\newcommand{\smj}{\sum_{j=1}^{m}}
\newcommand{\sumjm}{\sum_{j=1}^{m}}
\newcommand{\bdis}{\begin{displaymath}}
\newcommand{\edis}{\end{displaymath}}
\newcommand{\beq}{\begin{equation}}
\newcommand{\eeq}{\end{equation}}
\newcommand{\beqn}{\begin{eqnarray}}
\newcommand{\eeqn}{\end{eqnarray}}
\newcommand{\defeq}{\stackrel{\triangle}{=}}
\newcommand{\simleq}{\stackrel{<}{\sim}}
\newcommand{\sep}{\;\;\;\;\;\; ; \;\;\;\;\;\;}
\newcommand{\real}{\mbox{$ I \hskip -4.0pt R $}}
\newcommand{\complex}{\mbox{$ I \hskip -6.8pt C $}}
\newcommand{\integ}{\mbox{$ Z $}}
\newcommand{\realn}{\real ^{n}}
\newcommand{\sqrn}{\sqrt{n}}
\newcommand{\sqrtwo}{\sqrt{2}}
\newcommand{\prf}{{\bf Proof. }}

\newcommand{\onehlf}{\frac{1}{2}}
\newcommand{\thrhlf}{\frac{3}{2}}
\newcommand{\fivhlf}{\frac{5}{2}}
\newcommand{\onethd}{\frac{1}{3}}
\newcommand{\ndt}{\frac{n}{2}}

\newcommand{\lb}{\left ( }
\newcommand{\lcb}{\left \{ }
\newcommand{\lsb}{\left [ }
\newcommand{\labs}{\left | }
\newcommand{\rb}{\right ) }
\newcommand{\rcb}{\right \} }
\newcommand{\rsb}{\right ] }
\newcommand{\rabs}{\right | }
\newcommand{\lnm}{\left \| }
\newcommand{\rnm}{\right \| }
\newcommand{\lambdab}{\bar{\lambda}}
%
%
\newcommand{\xj}{x_{j}}
\newcommand{\xjb}{\bar{x}_{j}}
\newcommand{\xro}{x_{\resh}}
\newcommand{\xrob}{\bar{x}_{\resh}}
\newcommand{\xsig}{x_{\sigma}}
\newcommand{\xsigb}{\bar{x}_{\sigma}}
\newcommand{\xnmjb}{\bar{x}_{n-j+1}}
\newcommand{\xnmj}{x_{n-j+1}}
\newcommand{\aroj}{a_{\resh j}}
\newcommand{\arojb}{\bar{a}_{\resh j}}
\newcommand{\aroro}{a_{\resh \resh}}
\newcommand{\amuro}{a_{\mu \resh}}
\newcommand{\amumu}{a_{\mu \mu}}
\newcommand{\aii}{a_{ii}}
\newcommand{\aik}{a_{ik}}
\newcommand{\akj}{a_{kj}}
\newcommand{\atwoii}{a^{(2)}_{ii}}
\newcommand{\atwoij}{a^{(2)}_{ij}}
\newcommand{\ajj}{a_{jj}}
\newcommand{\aiib}{\bar{a}_{ii}}
\newcommand{\ajjb}{\bar{a}_{jj}}
\newcommand{\bii}{a_{jj}}
\newcommand{\biib}{\bar{a}_{jj}}
\newcommand{\aij}{a_{i,n-i+1}}
\newcommand{\akl}{a_{j,n-j+1}}
\newcommand{\aijb}{\bar{a}_{i,n-i+1}}
\newcommand{\aklb}{\bar{a}_{j,n-j+1}}
\newcommand{\bij}{a_{n-j+1,j}}
\newcommand{\arorob}{\bar{a}_{\resh \resh}}
\newcommand{\arosig}{a_{\resh \sigma}}
\newcommand{\arosigb}{\bar{a}_{\resh \sigma}}
\newcommand{\sumjrosig}{\sum_{\stackrel{j=1}{j\neq\resh,\sigma}}^{n}}
\newcommand{\summuro}{\sum_{\stackrel{j=1}{j\neq\mu,\resh}}^{n}}
\newcommand{\sumjnoti}{\sum_{\stackrel{j=1}{j\neq i}}^{n}}
\newcommand{\sumlnoti}{\sum_{\stackrel{\ell=1}{\ell \neq i}}^{n}}
\newcommand{\sumknoti}{\sum_{\stackrel{k=1}{k\neq i}}^{n}}
\newcommand{\sumknotij}{\sum_{\stackrel{k=1}{k\neq i,j}}^{n}}
\newcommand{\sumk}{\sum_{k=1}^{n}}
\newcommand{\snl}{\sum_{\ell=1}^{n}}
\newcommand{\sumji}{\sum_{\stackrel{j=1}{j\neq i, n-i+1}}^{n}}
\newcommand{\sumki}{\sum_{\stackrel{k=1}{k\neq i, n-i+1}}^{n}}
\newcommand{\sumkj}{\sum_{\stackrel{k=1}{k\neq j, n-j+1}}^{n}}
\newcommand{\sumjro}{\sum_{\stackrel{j=1}{j\neq\resh}}^{n}}
\newcommand{\rrosig}{R''_{\resh \sigma}}
\newcommand{\rro}{R'_{\resh}}
\newcommand{\gamror}{\Gamma_{\resh}^{R}(A)}
\newcommand{\gamir}{\Gamma_{i}^{R}(A)}
\newcommand{\gamctrr}{\Gamma_{\frac{n+1}{2}}^{R}(A)}
\newcommand{\gamctrc}{\Gamma_{\frac{n+1}{2}}^{C}(A)}
\newcommand{\gamroc}{\Gamma_{\resh}^{C}(A)}
\newcommand{\gamjc}{\Gamma_{j}^{C}(A)}
\newcommand{\lamror}{\Lambda_{\resh}^{R}(A)}
\newcommand{\lamir}{\Lambda_{i}^{R}(A)}
\newcommand{\lamirepsilon}{\Lambda_{i}^{R}(A_{\epsilon})}
\newcommand{\lamnir}{\Lambda_{n-i+1}^{R}(A)}
\newcommand{\lamjr}{\Lambda_{j}^{R}(A)}
\newcommand{\varphiij}{\Phi_{ij}^{R}(A)}
\newcommand{\delir}{\Delta_{i}^{R}(A)}
\newcommand{\vir}{V_{i}^{R}(A)}
\newcommand{\pamir}{\Pi_{i}^{R}(A)}
\newcommand{\xir}{\Xi_{i}^{R}(A)}
\newcommand{\lamjc}{\Lambda_{j}^{C}(A)}
\newcommand{\vjc}{V_{j}^{C}(A)}
\newcommand{\pamjc}{\Pi_{j}^{C}(A)}
\newcommand{\xjc}{\Xi_{j}^{C}(A)}
\newcommand{\lamroc}{\Lambda_{\resh}^{C}(A)}
\newcommand{\lamsigr}{\Lambda_{\sigma}^{R}(A)}
\newcommand{\lamsigc}{\Lambda_{\sigma}^{C}(A)}
\newcommand{\psii}{\Psi_{i}^{R}(A)}
\newcommand{\psiq}{\Psi_{q}}
\newcommand{\psiiepsilon}{\Psi_{i}^{R}(A_{\epsilon})}
\newcommand{\psiqepsilon}{\Psi_{q}(A_{\epsilon})}
\newcommand{\psiqc}{\Psi_{q}^{c}}
\newcommand{\psiqcepsilon}{\Psi_{q}^{c}(A_{\epsilon})}

\newcommand{\xmu}{x_{\mu}}
\newcommand{\xmub}{\bar{x}_{\mu}}
\newcommand{\xnu}{x_{\nu}}
\newcommand{\xnub}{\bar{x}_{\nu}}
\newcommand{\amuj}{a_{\mu j}}
\newcommand{\amujb}{\bar{a}_{\mu j}}
\newcommand{\amumub}{\bar{a}_{\mu \mu}}
\newcommand{\amunu}{a_{\mu \nu}}
\newcommand{\amunub}{\bar{a}_{\mu \nu}}
\newcommand{\sumjmunu}{\sum_{\stackrel{j=1}{j\neq\mu,\nu}}}
\newcommand{\rmunu}{R''_{\mu \nu}}
\newcommand{\rmu}{R'_{\mu}}

\newcommand{\Azero}{A_{0}}
\newcommand{\Aone}{A_{1}}
\newcommand{\Atwo}{A_{2}}
\newcommand{\Ath}{A_{3}}
\newcommand{\Afr}{A_{4}}
\newcommand{\Afv}{A_{5}}
\newcommand{\Asx}{A_{6}}
\newcommand{\Anmo}{A_{n-1}}
\newcommand{\Anmt}{A_{n-2}}
\newcommand{\Ajmo}{A_{j-1}}
\newcommand{\Ajpo}{A_{j+1}}
\newcommand{\Aj}{A_{j}}
\newcommand{\Andtmo}{A_{\frac{n}{2}-1}}
\newcommand{\Ammo}{A_{m-1}}
\newcommand{\Am}{A_{m}}
\newcommand{\Andt}{A_{\frac{n}{2}}}
\newcommand{\An}{A_{n}}

\newcommand{\azero}{a_{0}}
\newcommand{\aone}{a_{1}}
\newcommand{\atwo}{a_{2}}
\newcommand{\ath}{a_{3}}
\newcommand{\afr}{a_{4}}
\newcommand{\afv}{a_{5}}
\newcommand{\asx}{a_{6}}
\newcommand{\anmo}{a_{n-1}}
\newcommand{\anmt}{a_{n-2}}
\newcommand{\anmth}{a_{n-3}}
\newcommand{\ajmo}{a_{j-1}}
\newcommand{\ajpo}{a_{j+1}}
\newcommand{\aj}{a_{j}}
\newcommand{\andtmo}{a_{\frac{n}{2}-1}}
\newcommand{\andt}{a_{\frac{n}{2}}}
\newcommand{\an}{a_{n}}

\newcommand{\Bzero}{B_{0}}
\newcommand{\Bone}{B_{1}}
\newcommand{\Btwo}{B_{2}}
\newcommand{\Bth}{B_{3}}
\newcommand{\Bfr}{B_{4}}
\newcommand{\Bfv}{B_{5}}
\newcommand{\Bsx}{B_{6}}
\newcommand{\Bnmo}{B_{n-1}}
\newcommand{\Bndtwomo}{B_{n/2-1}}
\newcommand{\Bnmt}{B_{n-2}}
\newcommand{\Bjmo}{B_{j-1}}
\newcommand{\Bjpo}{B_{j+1}}
\newcommand{\Bj}{B_{j}}
\newcommand{\Bndtmo}{B_{\frac{n}{2}-1}}
\newcommand{\Bndt}{B_{\frac{n}{2}}}
\newcommand{\Bn}{B_{n}}

\newcommand{\vzero}{v_{0}}
\newcommand{\vone}{v_{1}}
\newcommand{\vtwo}{v_{2}}
\newcommand{\vj}{v_{j}}
\newcommand{\vnmth}{v_{n-3}}
\newcommand{\vnmtwo}{v_{n-2}}
\newcommand{\vnmo}{v_{n-1}}
\newcommand{\vn}{v_{n}}

\newcommand{\wzero}{w_{0}}
\newcommand{\wone}{w_{1}}
\newcommand{\wtwo}{w_{2}}
\newcommand{\wj}{w_{j}}
\newcommand{\wnmth}{w_{n-3}}
\newcommand{\wnmtwo}{w_{n-2}}
\newcommand{\wnmo}{w_{n-1}}
\newcommand{\wn}{w_{n}}

\newcommand{\tildet}{\tilde{t}}
\newcommand{\tildeu}{\tilde{u}}
\newcommand{\tildemu}{\tilde{\mu}}

\newcommand{\cii}{c_{ii}}
\newcommand{\cik}{c_{ik}}
\newcommand{\ckj}{c_{kj}}
\newcommand{\ctwoii}{c^{(2)}_{ii}}
\newcommand{\ctwoij}{c^{(2)}_{ij}}
\newcommand{\cjj}{c_{jj}}

\newcommand{\bik}{b_{ik}}
\newcommand{\bkj}{b_{kj}}
\newcommand{\btwoii}{b^{(2)}_{ii}}
\newcommand{\btwoij}{b^{(2)}_{ij}}
\newcommand{\bjj}{b_{jj}}

\newcommand{\abii}{(AB)_{ii}}
\newcommand{\abil}{(AB)_{i\ell}}

\newcommand{\bkl}{b_{k\ell}}
\newcommand{\btwoil}{b^{(2)}_{i\ell}}
\newcommand{\bll}{b_{\ell \ell}}

\newcommand{\matrixspace}{\;\;}
\newcommand{\ellone}{\ell_{1}}
\newcommand{\elltwo}{\ell_{2}}

\newcommand{\varphik}{\varphi_{k}}
\newcommand{\chik}{\chi_{k}}
\newcommand{\Phik}{\Phi_{k}}
\newcommand{\psik}{\psi_{k}}
\newcommand{\dr}{\beta^{n-k}}
\newcommand{\dn}{\beta^{-k}}
\newcommand{\betatwotilde}{\tilde{\rtwonk}}
\newcommand{\betathtilde}{\tilde{\ronenk}}
\newcommand{\betamink}{\beta^{-k}}

\newcommand{\rvarphik}{r(\varphi_{k})}    
\newcommand{\rrvarphik}{(r(\varphi_{k}))}    
\newcommand{\rchik}{r(\chi_{k})}    
\newcommand{\rrchik}{(r(\chi_{k}))}    
\newcommand{\svarphik}{s(\varphi_{k})}    
\newcommand{\ssvarphik}{(s(\varphi_{k}))}    
\newcommand{\schik}{s(\chi_{k})}    
\newcommand{\sschik}{(s(\chi_{k}))}    

\newcommand{\rpsik}{r(\psi_{k})}    
\newcommand{\rrpsik}{(r(\psi_{k}))} 
\newcommand{\spsik}{s(\psi_{k})}    
\newcommand{\sspsik}{(s(\psi_{k}))} 
\newcommand{\rpsinmk}{r(\psi_{n-k})}    
\newcommand{\rrpsinmk}{(r(\psi_{n-k}))} 
\newcommand{\spsinmk}{s(\psi_{n-k})}    
\newcommand{\sspsinmk}{(s(\psi_{n-k}))} 

\newcommand{\rvarphinmk}{r(\varphi_{n-k})}    
\newcommand{\rrvarphinmk}{(r(\varphi_{n-k}))}    
\newcommand{\rchinmk}{r(\chi_{n-k})}    
\newcommand{\rrchinmk}{(r(\chi_{n-k}))}    
\newcommand{\svarphinmk}{s(\varphi_{n-k})}    
\newcommand{\ssvarphinmk}{(s(\varphi_{n-k}))}    
\newcommand{\schinmk}{s(\chi_{n-k})}    
\newcommand{\sschinmk}{(s(\chi_{n-k}))}    
\newcommand{\stilde}{\tilde{s}}

\newcommand{\resh}{\rho}
\newcommand{\snk}{s}
\newcommand{\ronenk}{r_{1}}
\newcommand{\rtwonk}{r_{2}}
\newcommand{\obar}{\bar{O}}
\newcommand{\Lbar}{\bar{L}}
\newcommand{\dltaone}{\delta_{1}}
\newcommand{\dltatwo}{\delta_{2}}
\newcommand{\dltatld}{\tilde{\delta}}
\newcommand{\tauone}{\tau_{1}}
\newcommand{\tautwo}{\tau_{2}}
\newcommand{\xb}{\bar{x}}
\newcommand{\qone}{q_{1}}
\newcommand{\qtwo}{q_{2}}
\newcommand{\ub}{\bar{u}}
\newcommand{\cmm}{\complex^{m \times m}}
\newcommand{\opdsk}{\mathcal{O}}
\newcommand{\cldsk}{\bar{\mathcal{O}}}

\newcommand{\sqrtaone}{\sqrt{\a1}}
\newcommand{\sqrtatwo}{\sqrt{\atwo}}

\newcommand{\rhoone}{\rho_{1}}
\newcommand{\rhotwo}{\rho_{2}}
\newcommand{\xone}{x_{1}}
\newcommand{\xtwo}{x_{2}}
\newcommand{\xthr}{x_{3}}
\newcommand{\xfor}{x_{4}}
\newcommand{\xfiv}{x_{5}}

\newcommand{\lnorm}{\left \|}
\newcommand{\rnorm}{\right \|}
\newcommand{\lnrm}{\biggl | \biggl |}
\newcommand{\rnrm}{\biggr |\biggr |}

\newcommand{\recipp}{p^{\protect \#}}
\newcommand{\recipP}{P^{\protect \#}}

\newcommand{\varphione}{\varphi_{1}}
\newcommand{\varphitwo}{\varphi_{2}}

\newcommand{\psione}{\psi_{1}}
\newcommand{\psitwo}{\psi_{2}}

\newcommand{\absatwo}{|\atwo|}
\newcommand{\absqrtatwo}{\sqrt{|\atwo|}}
\newcommand{\absrhoone}{|\rhoone|}
\newcommand{\absrhotwo}{|\rhotwo|}
\newcommand{\sqrgamma}{\sqrt{\gamma}}
\newcommand{\aminb}{a-b}                 
\newcommand{\sqrtc}{\sqrt{c}}            
\newcommand{\sqrtabsc}{\sqrt{|c|}}            
\newcommand{\absalf}{|\alpha|}            
\newcommand{\rhoi}{\rho_{i}}                  
\newcommand{\sigmai}{\sigma_{i}}                  
\newcommand{\rhoj}{\rho_{j}}                  
\newcommand{\sigmaj}{\sigma_{j}}                  
\newcommand{\xij}{x_{ij}}                  
\newcommand{\yij}{y_{ij}}                  
\newcommand{\varphii}{\varphi_{i}}
\newcommand{\varphij}{\varphi_{j}}
\newcommand{\xstar}{x^{\ast}}
\newcommand{\uuone}{\mathcal{U}_{1}}
\newcommand{\uutwo}{\mathcal{U}_{2}}
\newcommand{\xkova}{\hat{x}}
\newcommand{\xkovalam}{\xkova_{max}(\ell-1)}

\begin{center}
\large
{\bf CAUCHY-LIKE AND PELLET-LIKE RESULTS FOR POLYNOMIALS}
\vskip 0.5cm
\normalsize
A. Melman \\
Department of Applied Mathematics \\
School of Engineering, Santa Clara University  \\
Santa Clara, CA 95053  \\
e-mail : amelman@scu.edu \\
\vskip 0.5cm
\end{center}

\begin{abstract}
We obtain several Cauchy-like and Pellet-like results for the zeros of a general complex polynomial by considering similarity transformations 
of the squared companion matrix and the reformulation of the zeros of a scalar polynomial as the eigenvalues of a polynomial eigenvalue problem.
\vskip 0.15cm
{\bf Key words :} matrix polynomial, companion matrix, Cauchy, Pellet
\vskip 0.15cm
{\bf AMS(MOS) subject classification :} 12D10, 15A18, 30C15
\end{abstract}

%
%
                                                            \section{Introduction}           
                                                             \label{introduction}                                             
%
%

Two standard results for the localization of all or some of the zeros of a polynomial, due to Cauchy and Pellet, respectively, are given by
%
%
\begin{theorem}
{\bf (Cauchy's theorem - original scalar version)} (\cite{Cauchy}, \cite[Th.(27,1), p.122 and Exercise 1, p.126]{Marden})
\label{CauchyTheorem}             
All the zeros of the polynomial $p(z)= z^{n} + a_{n-1} z^{n-1} + \dots + a_{1} z + a_{0}$ with complex coefficients, $n \geq 2$,
lie in $|z| \leq s$, where $s$ is the unique positive solution of 
\bdis
x^{n} - |a_{n-1}| x^{n-1} - \dots  - |\aone| x - |\azero| = 0  \; . 
\edis
\end{theorem}

%
%
\begin{theorem}
{\bf (Pellet's theorem - original scalar version)} (\cite{Pellet}, \cite[Th.(28,1), p.128]{Marden})
\label{PelletTheorem}
Given the polynomial $p(z)= z^{n} + a_{n-1} z^{n-1} + \dots + a_{1} z + a_{0}$ with complex coefficients, 
$a_{k} \neq 0$, and $n \geq 2$, let $1 \leq k \leq n-1$, and let the polynomial 
\bdis  
x^{n} + |a_{n-1}| x^{n-1} + \dots + |a_{k+1}| x^{k+1} - |a_{k}| x^{k} + |a_{k-1}|x^{k-1} + \dots  + |a_{0}| 
\edis   
have two distinct positive roots $x_{1}$ and $x_{2}$ with $x_{1} < x_{2}$.
Then $p$ has exactly $k$ zeros in or on the circle $|z| = x_{1}$ and no zeros in the annular ring 
$x_{1} < |z| < x_{2}$.
\end{theorem}

Theorem~\ref{CauchyTheorem} provides an upper bound on the moduli of the zeros, whereas Theorem~\ref{PelletTheorem} sometimes allows zeros to be separated
into two different groups, according to the magnitude of their moduli. However, the latter is very sensitive to the magnitude of the coefficients and for the theorem
to be applicable, one or more coefficients typically have to be much larger than the others. 
 
The inequalities for the moduli of the zeros in these theorems are sharp in the sense that there exist polynomials for which they hold as equalities. 
Our aim is nevertheless to improve Theorem~\ref{CauchyTheorem} and a few special cases of Theorem~\ref{PelletTheorem}, 
resulting in a number of results of a similar nature, i.e., also involving the solution of one or two real equations. We immediately point out that the solution 
of such equations requires a negligible computational effort compared to the computation of the actual zeros 
(see, e.g., \cite{MelmanImplement}, \cite{Rump}), and we will not dwell on it.
Theorems~\ref{CauchyTheorem} and~\ref{PelletTheorem} have many applications and are often used to find good starting points for iterative methods that 
compute some or all of the zeros.

There are several ways to derive these and many other results related to polynomial zeros, one of which is to use linear algebra arguments.
Although not necessarily producing the shortest proofs, this provides a transparent and often elegant treatment of such results.
On the other hand, a linear algebra approach does not generally seem to lead to results that cannot also be obtained by purely
algebraic manipulation or applications of complex analysis, an observation also made in~\cite[p.263]{RS}.
Here, in contrast, we will use linear algebra tools to derive results that, it appears, cannot easily be obtained otherwise.      

Before going into more detail, we recall that the zeros of the complex monic scalar polynomial $p(z) = z^{n} + \anmo z^{n-1} + \dots + \azero$ 
are the eigenvalues of the $n \times n$ \emph{companion matrix} $C(p)$, defined by
\beq
\label{compmatdef}
C(p) =
\begin{pmatrix}
0 &        &       &   & -\azero     \\
1 &        &       &   & -\aone     \\
  & \ddots &       &   & \vdots \\
  &        &       & 1 & -\anmo    \\
\end{pmatrix}
\; .
\eeq    
Blank spaces in the matrices indicate zero elements. Thus, locating the eigenvalues of $C(p)$ is equivalent to locating the zeros of $p$. 
We will make frequent use of Gershgorin's theorem, which provides inclusion regions for the eigenvalues of a matrix. It is stated next.
%
%
\begin{theorem}
\label{GershgorinTheorem}
{\bf (Gershgorin's theorem)} (\cite{G}, \cite[Section 6.1]{HJ})
All the eigenvalues of the $n \times n$ complex matrix $A$ with elements $a_{ij}$ and deleted row sums $R'_{i}(A) = \sum_{\stackrel{j=1}{j\neq i}}^{n} |a_{ij} |$
are located in the union of $n$ discs $\bigcup_{i=1}^{n} \obar \lb a_{ii} ; R'_{i}(A) \rb $.
If $k$ discs are disjoint from the other discs, then their union contains exactly $k$ eigenvalues.
\end{theorem}
Because the eigenvalues of a matrix $A$ are the same as those of its transpose $A^{T}$, an analogous version is obtained by interchanging rows and columns. 
We refer to these as the row and column versions of the theorem and to the eigenvalue inclusion regions as the Gershgorin row and column sets, respectively. 
A good in-depth exposition of this theorem and many related theorems can be found in~\cite{VargaGershBook}.
As an illustration, let us apply the column version to $C(p)$: all the zeros of the polynomial $p$ 
can be found in the union of the disc, centered at the origin with radius one and the dick, centered at $-\anmo$, whose radius is the $n$th deleted column sum. 
It is often useful to apply an appropriate similarity transformation to the matrix, which does not change the eigenvalues, although it does affect the Gershgorin set.

For example, Theorem~\ref{CauchyTheorem} can be obtained by applying Gershgorin's theorem to a specific diagonal similarity transformation 
of the companion matrix (\cite{Bell}, \cite{MelmanTwin}, and \cite[Theorem 8.6.3]{RS}), and the same is true for related results (\cite{MelmanTwin}). 
We propose to improve the aforementioned results and derive additional ones by considering the square of the companion matrix, the eigenvalues of which are 
the squares of the zeros of $p$. It is given by
\beq     
\label{sqcompmatdef}
C^{2}(p) = 
\begin{pmatrix}
0 & 0  & \dots & 0 & -\azero  & \anmo \azero\\
0 & 0  & \dots & 0 & -\aone  & \anmo \aone-\azero\\
1 & 0  & \dots & 0 & -\atwo  & \anmo \atwo-\aone\\
\vdots & \vdots  & \vdots & \vdots & \vdots & \vdots\\
0 & 0  & \dots & 1 & -\anmo & \anmo^{2}-\anmt\\
\end{pmatrix} \; .
\eeq
The idea of obtaining additional inclusion regions for the eigenvalues of a general matrix by squaring it is certainly not new, but it    
does not typically lead to an improvement (for some examples see, e.g., \cite{MelmanCassiniToeplitz}). However, there are good reasons 
to use $C^{2}(p)$ instead of $C(p)$. First of all, it can be shown (\cite[Theorem 2.1]{MelmanCassiniToeplitz}) that squaring does lead 
to smaller inclusion regions when the matrix has a zero diagonal, and that is almost the case for a companion matrix: only one diagonal element is not
necessarily zero.
Secondly, of equal importance is the more complicated structure of the squared companion matrix while still keeping it relatively simple (only two columns). 
The simplicity of the companion matrix 
is an advantage when computing its eigenvalues, but it also means that linear algebra tools have less room to maneuver when trying to extract 
information about the location of the eigenvalues without actually computing them. As we will see, squaring expands the 
range of useful similarity transformations significantly, while also suggesting a natural and convenient reformulation of the zeros of a scalar polynomial
as the eigenvalues of a matrix polynomial, leading to further advantages.

Matrix polynomials are encountered when a nonzero complex vector $v$ and a complex number $z$ are sought such 
that $P(z)v=0$, with           
\bdis
P(z) = A_{n} z^{n} + A_{n-1} z^{n-1} + \dots + A_{0} \; ,
\edis 
and the coefficients $A_{j}$ are $m \times m$ complex matrices. If $A_{n}$ is singular then there are infinite eigenvalues and if $A_{0}$ is singular 
then zero is an eigenvalue. There are $nm$ eigenvalues, including possibly infinite ones. The finite eigenvalues are the solutions of $det{P(z)}=0$.
If $P$ is a monic matrix polynomial, i.e., $\An = I$, then its eigenvalues are the eigenvalues of the $nm \times nm$ \emph{block companion matrix} $C(P)$, 
defined by
\bdis
C(P) =
\begin{pmatrix}
0 &        &       &   & -\Azero     \\
I &        &       &   & -\Aone     \\
  & \ddots &       &   & \vdots \\
  &        &       & I & -\Anmo    \\
\end{pmatrix}
\; . 
\edis 
Since the size of $I$ will usually be clear from the context, we omit it from the notation.
Theorem~\ref{CauchyTheorem} and Theorem~\ref{PelletTheorem} have analogs for matrix polynomials, and they are stated next. 
The matrix norms are assumed to be subordinate (induced by a vector norm).

%
%
\begin{theorem}
{\bf (Cauchy's theorem - matrix version)} (\cite{BiniNoferiniSharify}, \cite{HighamTisseur}, \cite{Melman_MatPol})
\label{MatrixCauchyTheorem}             
All eigenvalues of the matrix polynomial $P(z) = Iz^{n} + A_{n-1}z^{n-1} + \dots + A_{1}z + A_{0}$, where $n \geq 2$ and $A_{j} \in \cmm$, for $j=0,\dots,n$, 
lie in $|z| \leq s$, where $s$ is the unique positive solution of
\bdis
x^{n} - ||A_{n-1}||x^{n-1} - \dots - ||A_{1}||x - ||A_{0}|| = 0 \; . 
\edis
\end{theorem}

%
%
\begin{theorem}
\label{MatrixPelletTheorem}             
{\bf (Pellet's theorem - matrix version.)} (\cite{BiniNoferiniSharify}, \cite{Melman_MatPol})
Let 
\bdis  
P(z) = I z^{n} + A_{n-1}z^{n-1} + \dots + A_{1}z + A_{0}  \nonumber 
\edis
be a matrix polynomial with $n \geq 2$, $A_{j} \in \cmm$ for $j=0,\dots,n-1$, and $A_{0} \neq 0$.      
Let $A_{k}$ be invertible for some $k$ with $1 \leq k \leq n-1$, and let the polynomial 
\bdis
x^{n} + ||A_{n-1}||x^{n-1} + \dots + ||A_{k+1}||x^{k+1} - ||A^{-1}_{k}||^{-1}x^{k} 
+ ||A_{k-1}||x^{k-1} + \dots + ||A_{1}||x + ||A_{0}||   
\edis
have two distinct positive roots $x_{1}$ and $x_{2}$ with $x_{1} < x_{2}$.
Then $\det(P)$ has exactly $km$ zeros in or on the disc $|z|=x_{1}$ and no zeros in the annular ring 
$x_{1} < |z| < x_{2}$.
\end{theorem}

Finally, we mention a theorem, which we call the Block Gershgorin theorem, due to D.G.~Feingold and R.S.~Varga (\cite{FV}). 
%
%
\begin{theorem}
\label{BlockGershgorinTheorem}
{\bf (Block Gershgorin theorem)} (\cite[Theorems 2 and 4]{FV})
Let $A$ be any $n \times n$ matrix with complex entries, which is partitioned in the following manner:
\bdis
\begin{pmatrix}
A_{1,1} & A_{1,2} & \dots  & A_{1,N} \\
A_{2,1} & A_{2,2} & \dots  & A_{2,N} \\
\vdots  & \vdots  & \vdots & \vdots  \\
A_{N,1} & A_{N,2} & \dots  & A_{N,N} \\
\end{pmatrix}
\; ,
\edis
where the diagonal submatrices $A_{i,i}$ are square of order $n_{i}$, $1 \leq i \leq N$, and let $I_{i}$ be the $n_{i} \times n_{i}$ identity matrix.
Then each eigenvalue $\lambda$ of $A$ lies in a Gershgorin set $G_{j}$ for at least one $j$, $1 \leq j \leq N$, where $G_{j}$ is the set of all complex 
numbers $z$ such that
\bdis
\left \| \bigl ( A_{j,j} - z I_{i} \bigr )^{-1} \right \|^{-1} \leq \sum_{\stackrel{k=1}{k \neq j}}^{N} ||A_{j,k}|| \; .
\edis
If the union $H=\bigcup_{j=1}^{k}G_{p_{j}}$, $1 \leq p_{j} \leq N$, of $k$ Gershgorin sets $G_{j}$ is disjoint from the remaining $N-k$ Gershgorin sets, then
$H$ contains precisely $\sum_{j=1}^{k} n_{p_{j}}$ eigenvalues of $A$.  
\end{theorem}
The norms in this theorem are subordinate and the theorem has, just like Gershgorin's theorem, also a block-column form. The sets $G_{j}$ are difficult to 
compute in general, but each is a union of $n_{j}$ discs when the Euclidean norm (2-norm) is used and the diagonal blocks are 
normal matrices (\cite[Theorem 5]{FV}). When the diagonal blocks are diagonal matrices, the sets are discs for any subordinate norm.

Lower bounds on the moduli of polynomial zeros can be obtained by applying the appropriate aforementioned theorems to the reverse polynomial $p^{\protect\#}$
of $p$, defined by $p^{\protect \#}(z) = z^{n}p(1/z)/\azero$, whose zeros are the reciprocals of those of $p$.

In the next section, we consider several similarity transformations of $C^{2}(p)$ that will be used in Section~3 to derive Cauchy-like results as in
Theorems~\ref{CauchyTheorem} and~\ref{MatrixCauchyTheorem}, and in Section~4 to derive similar results to those of Theorems~\ref{PelletTheorem} 
and~\ref{MatrixPelletTheorem}. We present numerical results in Section~5 to illustrate the (sometimes drastic) improvements we were able to obtain.
An effort was made to make statements of theorems self-contained with the inevitable repetition of some definitions.

%
%
                                                            \section{Preliminaries}    
                                                             \label{Preliminaries}                                            
%
%

Throughout, we denote by $O \lb a;r \rb$ and $\obar \lb a;r \rb$ the open and closed discs, respectively, centered at $a$ with radius $r$. 
Consider the polynomial $p(z)= z^{\ell} + b_{\ell-1} z^{\ell-1} + \dots + b_{1} z + b_{0}$. If it is not of even degree, then 
we consider $zp(z)$, which has a zero constant term. The latter has no effect on theoretical results, nor does it affect, as we will see, numerical results. 
To cover both cases at once, we define $m=\ceil*{\frac{\ell}{2}}$, $n=2m$, and the complex numbers $a_{j}$ as follows:
\begin{eqnarray*}
& & a_{j} = b_{j} \;\; \text{($\ell$ is even and $j=0,\dots,n-1$),}  \\
& & a_{0} = 0 \; \text{and} \; a_{j} = b_{j-1} \;\; \text{($\ell$ is odd and $j=1,\dots,n-1$).}  
\end{eqnarray*}
This means that, when $\ell$ is odd, we have $n=\ell+1$ and $m=(\ell + 1)/2$. When $\ell$ is even, $n=\ell$ and $m=\ell/2$.
From here on, we consider the polynomial $p(z)= z^{n} + a_{n-1} z^{n-1} + \dots + a_{1} z + a_{0}$ of even degree $n$. 
From~(\ref{sqcompmatdef}), the square of the companion matrix of $p$ can be written as
\beq      
\nonumber
C^{2}(p) = 
\begin{pmatrix}
0 & 0  & \dots & 0 & -\azero  & \anmo \azero\\
0 & 0  & \dots & 0 & -\aone  & \anmo \aone-\azero\\
1 & 0  & \dots & 0 & -\atwo  & \anmo \atwo-\aone\\
\vdots & \vdots  & \vdots & \vdots & \vdots & \vdots\\
0 & 0  & \dots & 1 & -\anmo & \anmo^{2}-\anmt\\
\end{pmatrix} 
=
\begin{pmatrix}
0 &        &       &   & -\Azero     \\
I &        &       &   & -\Aone     \\
  & \ddots &       &   & \vdots \\
  &        &       & I & -\Ammo    \\
\end{pmatrix}
\; ,
\eeq
where
\bdis 
A_{0} = \begin{pmatrix} -\azero & \anmo \azero \\ -\aone & \anmo \aone -\azero \end{pmatrix} \; , 
\;\; \text{and} \;\;
A_{j} = \begin{pmatrix} -a_{2j} & \anmo a_{2j} -a_{2j-1} \\ -a_{2j+1} & \anmo a_{2j+1} -a_{2j} \\ \end{pmatrix} 
\;\; \text{for $j=1,...,m-1$.}
\edis

By Schur's theorem, there exists a unitary matrix $U$ that triangularizes the matrix 
\bdis
\Ammo = \begin{pmatrix} -a_{n-2} & \anmo a_{n-2} -a_{n-3} \\ -a_{n-1} & \anmo^{2} -a_{n-2} \\ \end{pmatrix} \; , 
\edis
i.e., 
\bdis
U^{*}\Ammo U = \begin{pmatrix} \alpha & \gamma \\ 0 & \beta \\ \end{pmatrix} \; ,
\edis
where $\alpha, \beta, \gamma \in \complex$, $\alpha$ and $\beta$ are the eigenvalues of $\Ammo$, and $U^{*}=U^{-1}$ is the Hermitian conjugate of $U$.
If $\Ammo$ is diagonalizable, then there exists a nonsingular matrix $M$ for which
\bdis
M^{-1}\Ammo M = \begin{pmatrix} \alpha & 0 \\ 0 & \beta \\ \end{pmatrix} \; .
\edis
In what follows, the matrix $S$ is defined as a matrix that either triangularizes $\Ammo$ or diagonalizes it if that is possible. 

We now consider three similarity transformations of $C^{2}(p)$ together with
their Gershgorin column sets, that will form the building blocks for the Cauchy-like and Pellet-like results of the following sections. 

Let $Q$ be an $n \times n$ block diagonal matrix with $m=n/2$ identical diagonal blocks equal to $S$. Then we define $C_{S}^{2}(p) = Q^{-1}C^{2}(p)Q$, and,
with~(\ref{sqcompmatdef}), obtain
\bdis
C^{2}_{S}(p)
=
\begin{pmatrix}
0 &        &       &   & S^{-1} \Azero S   \\
I &        &       &   & S^{-1} \Aone S    \\
  & \ddots &       &   & \vdots            \\
  &        &       & I & S^{-1} \Ammo S    \\
\end{pmatrix}
=
\begin{pmatrix}
0   & 0    &    &        &     &     &           &  \vzero   & \wzero            \\
0   & 0    &    &        &     &     &           &  \vone    & \wone             \\
1   & 0    &    &        &     &     &           &  \vtwo    & \wtwo             \\
0   & 1    &    &        &     &     &           &  v_{3}    & w_{3}             \\
    &      &    & \ddots &     &     &           &  \vdots   & \vdots            \\
0   & 0    &    &        &     & 1   & 0         &  \alpha   & \gamma            \\
0   & 0    &    &        &     & 0   & 1         &  0        & \beta             \\
\end{pmatrix}
\; ,
\edis  
where the vectors $v,w \in \complex^{n-2}$ are defined by                   
\bdis  
\begin{pmatrix} v_{2j} & w_{2j} \\ v_{2j+1} & w_{2j+1} \\ \end{pmatrix} = S^{-1} \Aj S  
\;\; \text{for $j=0,...,m-2$} \; ,
\edis           
and $\alpha , \beta , \gamma \in \complex$ with $\gamma = 0$ if $\Ammo$ is diagonalizable.
The triangularization (or diagonalization) of the lower right-hand block of $C^{2}(p)$ will facilitate the application of Gershgorin's theorem.
 
To $C^{2}_{S}(p)$ we apply two different similarity transformations. First, let $D_{x}$ be the diagonal matrix with diagonal ($x^{n},x^{n-1},\dots,x$), and define
\beq
\label{Fdef}
F_{x}(p) = D_{x}^{-1} C^{2}_{S}(p) D_{x}
=
\begin{pmatrix}
0       & 0      &    &        &     &        &           &  \vzero/x^{n-2}   & \wzero/x^{n-1}            \\
0       & 0      &    &        &     &        &           &  \vone/x^{n-3}    & \wone/x^{n-2}             \\
x^{2}   & 0      &    &        &     &        &           &  \vtwo/x^{n-4}    & \wtwo/x^{n-3}             \\
0       & x^{2}  &    &        &     &        &           &  v_{3}/x^{n-5}    & w_{3}/x^{n-4}             \\
        &        &    & \ddots &     &        &           &  \vdots   & \vdots                            \\
0       & 0      &    &        &     & x^{2}  & 0         &  \alpha   & \gamma/x                          \\
0       & 0      &    &        &     & 0      & x^{2}     &  0        & \beta                             \\
\end{pmatrix}
\; .
\eeq
Then for any $x>0$, the Gershgorin column set of $F_{x}(p)$ is the union of the three discs $O_{1}(x) \equiv \obar(0;x^{2})$,
$O_{2}(x) \equiv \obar(\alpha;\rho_{1}(x))$, and $O_{3}(x) \equiv \obar(\beta;\rho_{2}(x))$, where
\begin{eqnarray}
& & \rho_{1}(x) = \dfrac{|v_{n-3} |}{x} + \dfrac{|v_{n-4} |}{x^{2}} + \dots + \dfrac{|\vone |}{x^{n-3}} + \dfrac{|\vzero |}{x^{n-2}} \; ,  \label{rho1}  \\ 
\text{and} & &  \nonumber  \\
& & \rho_{2}(x) = \dfrac{|\gamma|}{x} + \dfrac{|\wnmth |}{x^{2}} + \dfrac{|w_{n-4} |}{x^{3}} + \dots + \dfrac{|\wone |}{x^{n-2}} 
+ \dfrac{|\wzero |}{x^{n-1}}  \; . \label{rho2} 
\end{eqnarray}
As $x$ varies, the discs expand and contract, and we will use this flexibility later to obtain convenient configurations of the discs.

Secondly, let $\Delta_{x}$ be the block diagonal matrix with diagonal blocks ($x^{m}I,x^{m-1}I,\dots,xI$), where $I$ is the $2 \times 2$ identity matrix
and $x>0$. We define
\beq     
\label{Phidef}
\Phi_{x}(p) = \Delta_{x}^{-1} C^{2}_{S}(p) \Delta_{x}
=
\begin{pmatrix}
0       & 0      &    &        &     &        &           &  \vzero/x^{m-1}   & \wzero/x^{m-1}            \\
0       & 0      &    &        &     &        &           &  \vone/x^{m-1}    & \wone/x^{m-1}             \\
x       & 0      &    &        &     &        &           &  \vtwo/x^{m-2}    & \wtwo/x^{m-2}             \\
0       & x      &    &        &     &        &           &  v_{3}/x^{m-2}    & w_{3}/x^{m-2}             \\
        &        &    & \ddots &     &        &           &  \vdots   & \vdots            \\
0       & 0      &    &        &     & x      & 0         &  \alpha   & \gamma            \\
0       & 0      &    &        &     & 0      & x         &  0        & \beta             \\
\end{pmatrix}
\; .
\eeq
The Gershgorin column set of $\Phi_{x}(p)$ is the union of the three discs $O'_{1}(x) \equiv \obar(0;x)$,
$O'_{2}(x) \equiv \obar(\alpha;\sigma_{1}(x))$, and $O'_{3}(x) \equiv \obar(\beta;\sigma_{2}(x))$, where
\begin{eqnarray}
& & \sigma_{1}(x) = \dfrac{|\vnmth| + |v_{n-4} |}{x} + \dfrac{|v_{n-5}| + |v_{n-6}|}{x^{2}} + \dots + \dfrac{|\vone | + |\vzero |}{x^{m-1}} \; , \label{sigma1} \\ 
\text{and} & &  \nonumber  \\
& & \sigma_{2}(x) = |\gamma| + \dfrac{|\wnmth | + |w_{n-4} |}{x} + \dfrac{|w_{n-5}| + |w_{n-6}|}{x^{2}} + \dots + \dfrac{|\wone | 
+ |\wzero |}{x^{m-1}} \; .  \label{sigma2} 
\end{eqnarray}

We define
\bdis
T_{m-1} = \begin{pmatrix} \alpha & \gamma \\ 0 & \beta \\ \end{pmatrix} 
\;\; , \;\; 
T_{j} = S^{-1} \Aj S   \;\; \text{for $j=0,...,m-2$,}
\edis
and, for any subordinate matrix norm,
\bdis
\tau(x) = \dfrac{||T_{m-2}||}{x} + \dfrac{||T_{m-3}||}{x^{2}} + \dots + \dfrac{||T_{0}||}{x^{m-1}} \; .
\edis
If $\Ammo$ is diagonalizable, we choose the matrix $S$ such that $\gamma = 0$ and  
the Block Gershgorin column set of $\Phi_{x}(p)$ is then given by $G_{1} \cup G_{2}$, 
where the sets $G_{j}$ are as defined in the statement of Therorem~\ref{BlockGershgorinTheorem}. It is straightforward to show that 
$G_{1} = \obar(0;x)$ and that 
\bdis
G_{2} = \left \{ z : \left \| \left ( \begin{pmatrix} \alpha & 0 \\ 0 & \beta \\ \end{pmatrix}  
- z \begin{pmatrix} 1 & 0 \\ 0 & 1 \\ \end{pmatrix} \right )^{-1} \right \|^{-1} \leq \tau(x) \right \}
= \left \{ z : \max{|z-\alpha|,|z-\beta|} \leq \tau(x) \right \} \; .
\edis
This means that $G_{2} = \obar(\alpha;\tau(x)) \cup \obar(\beta;\tau(x))$, the union of two discs
centered at $\alpha$ and $\beta$, respectively, with identical radii $\tau(x)$.
For this theorem, it is important for $\Ammo$ to be diagonalizable since the inclusion region would otherwise become too complicated to be useful.

\noindent {\bf Remarks.} 
\begin{itemize}
\item
One observes that $C^{2}(p)$ is the block companion matrix of the matrix polynomial 
$P(z) = I z^{m} + \Ammo z^{n-1} + \dots + \Aone z + \Azero$, i.e., the squares of the zeros of $p$ are also the eigenvalues of $P$.
It follows that they are also the eigenvalues of 
\beq
\label{PSdef}
P_{S}(z) = I z^{m} + T_{m-1} z^{n-1} + \dots + T_{1} z + T_{0} \; .
\eeq
\item
Since it simplifies all the equations we will encounter, we will diagonalize $\Ammo$ whenever this is possible.
It is a straightforward exercise to determine when $\Ammo$ is not diagonalizable in terms of the coefficients of $p$, since it can be 
written as
\bdis
A_{m-1} = \begin{pmatrix} -\anmt & \anmo \anmt - \anmth \\ -\anmo & \anmo^{2} -\anmt \\ \end{pmatrix} 
=
-\anmt I +  \begin{pmatrix}  0 & \anmo \anmt - \anmth \\ -\anmo & \anmo^{2}  \\ \end{pmatrix} 
\; .
\edis
Excluding $\anmo = \anmth = 0$ (which makes $\Ammo$ diagonal), $\Ammo$ cannot be diagonalized if it has a double eigenvalue.
Since the characteristic polynomial of $\Ammo + \anmt I$ is given by 
\bdis
\lambda^{2} - \anmo^{2} \lambda + \anmo(\anmo\anmt-\anmth) \; ,
\edis
the matrix $\Ammo$ is not diagonalizable when
\bdis
\anmo(\anmo^{3} - 4\anmo\anmt + 4 \anmth) = 0 \; ,
\edis
unless $\anmo = \anmth = 0$.
\item
As always when deriving bounds, an eye should be kept on the computational cost this entails, as this cost should remain well below the cost to 
compute the zeros exactly. This is certainly the case here, as the similarity transformations and the matrix norms, involving at most $2 \times 2$ blocks, 
require a total $\mathcal{O}(n)$ arithmetic operations, as do the solutions of the various real polynomial equations. The latter should require few iterations
when a properly adapted iterative method is used.
\end{itemize}
%
%
                                                            \section{Cauchy-like results}
                                                             \label{CauchyResults}                                            
%
%
We now present a theorem containing two Cauchy-like results for the moduli of a polynomial's zeros, using similarity transformations 
of the squared companion matrix.

%
%
\begin{theorem}
\label{TGC}                       
For a polynomial $p(z)= z^{\ell} + b_{\ell-1} z^{\ell-1} + \dots + b_{1} z + b_{0}$ with complex coefficients, $\ell \geq 3$, 
and zeros $\{z_{j} \}_{j=1}^{\ell}$, define $m = \ceil*{\frac{\ell}{2}}$, $n=2m$, and the complex numbers $a_{j}$ as follows:
\begin{eqnarray*}
& & a_{j} = b_{j} \;\; \text{($\ell$ is even and $j=0,\dots,n-1$),}  \\
& & a_{0} = 0 \; \text{and} \; a_{j} = b_{j-1} \;\; \text{($\ell$ is odd and $j=1,\dots,n-1$).}  
\end{eqnarray*}
Furthermore, define 
\bdis
A_{0} = \begin{pmatrix} -\azero & \anmo \azero \\ -\aone & \anmo \aone -\azero \end{pmatrix} \; , 
\;\; \text{and} \;\;
A_{j} = \begin{pmatrix} -a_{2j} & \anmo a_{2j} -a_{2j-1} \\ -a_{2j+1} & \anmo\ajpo -a_{2j} \\ \end{pmatrix} 
\;\; \text{for $j=1,...,m-1$,}
\edis
and let $S$ be a nonsingular matrix such that 
$            
S^{-1}\Ammo S = \begin{pmatrix} \alpha & \gamma \\ 0 & \beta \\ \end{pmatrix} \; ,
$           
where $\alpha, \beta, \gamma \in \complex$. Define the vectors $v,w \in \complex^{n-2}$ as                      
\bdis
\begin{pmatrix} v_{2j} & w_{2j} \\ v_{2j+1} & w_{2j+1} \\ \end{pmatrix} = S^{-1} \Aj S  
\;\; \text{for $j=0,...,m-2$} \; .
\edis

Then the following holds.
\newline {\bf (a)} 
$|z_{j}| \leq \max \left \{ r_{1} , r_{2} \right \}$,
where $r_{1}$ and $r_{2}$ are the unique positive solutions of $\psi_{1}(x) = 0$ and $\psi_{2}(x) = 0$, respectively, given by     
\begin{eqnarray*}
& & \psi_{1}(x) = x^{n} - |\alpha|x^{n-2} - |\vnmth |x^{n-3} - \dots - |\vone | x - |\vzero | \; ,      \label{psi1}   \\    
\text{and} & &  \\
& & \psi_{2}(x) = x^{n+1} - |\beta|x^{n-1} - |\gamma| x^{n-2} - |\wnmth |x^{n-3} - \dots - |\wone | x - |\wzero |   \; .   \label{psi2}   \\
\end{eqnarray*}
\noindent {\bf (b)} 
$|z_{j}| \leq \bigl ( \max \left \{ s_{1} , s_{2} \right \} \bigr )^{1/2}$,
where $s_{1}$ and $s_{2}$ are the unique positive solutions of $\varphi_{1}(x) = 0$ and $\varphi_{2}(x) = 0$, respectively, given by    
\begin{eqnarray*}
& & 
\hskip -0.75cm \varphi_{1} (x) = x^{m} - |\alpha|x^{m-1} - \big ( |\vnmth | + |v_{n-4} | \bigr ) x^{m-2}  
- \dots - \bigl ( |v_{3}| + |\vtwo| \bigr ) x - \dots - \lb |\vone | + |\vzero | \rb  \; ,   \\
& & \hskip -1.15cm \text{and}   \\
& & 
\hskip -0.75cm \varphi_{2}(x) = x^{m} - (|\beta|+|\gamma|)x^{m-1} - \big ( |\wnmth | + |w_{n-4} | \bigr ) x^{m-2}  
- \dots - \bigl ( |w_{3}| + |\wtwo| \bigr ) x - \dots - \lb |\wone | + |\wzero | \rb  \; .
\end{eqnarray*}
\end{theorem}
\noindent \prf
\newline {\bf (a)}
First assume that $n$ is even, in which case $m=\ell/2$, $n=\ell$, and $a_{j}=b_{j}$. Define $C^{2}_{S}(p)$ as in~(\ref{sqcompmatdef}), and for $x > 0$ let
$D_{x}$ be the diagonal matrix with diagonal ($x^{n},x^{n-1},\dots,x$). Define $F_{x}(p) = D_{x}^{-1} C^{2}_{S}(p) D_{x}$, 
so that, with~(\ref{Fdef}), the
Gershgorin column set of $F_{x}(p)$ is the union of the three discs $O_{1}(x) \equiv \obar(0;x^{2})$,
$O_{2}(x) \equiv \obar(\alpha;\rho_{1}(x))$, and $O_{3}(x) \equiv \obar(\beta;\rho_{2}(x))$, where $\rho_{1}$ and $\rho_{2}$ are defined by~(\ref{rho1})
and~(\ref{rho2}), respectively.
As $x$ increases, $O_{1}(x)$ expands, while $O_{2}(x)$ and $O_{3}(x)$ contract. When $x^{2} = |\alpha| + \rho_{1}(x)$, $O_{1}(x)$ and $O_{2}(x)$ will be 
tangent to one another and $O_{2}(x) \subseteq O_{1}(x)$. This happens when $x=r_{1}$, where $r_{1}$ is the unique positive solution of 
\bdis
x^{n} - |\alpha|x^{n-2} - |\vnmth |x^{n-3} - \dots - |\vone | x - |\vzero |  =  0   \; .
\edis
If $O_{3}(r_{1}) \subseteq O_{1}(r_{1})$, then the Gershgorin set is $O_{1}(r_{1})$. If this is not the case, we let $x$ increase further, until
$O_{1}(x)$ becomes tangent to $O_{3}(x)$ and $O_{3}(x) \subseteq O_{1}(x)$. This occurs when $x^{2} = |\beta| + \rho_{2}(x)$, which is when $x=r_{2}$, where 
$r_{2}$ is the unique positive solution of 
\bdis
x^{n+1} - |\beta|x^{n-1} - |\gamma| x^{n-2} - |\wnmth |x^{n-3} - \dots - |\wone | x - |\wzero |  =  0   \; .
\edis
The Gershgorin set of $F_{x}(p)$ is then equal to $O_{1}(r_{2})$. The case where $O_{1}(x)$ first becomes tangent to $O_{3}(x)$ is analogous.
We conclude that the Gershgorin set of $F_{x}(p)$ is given by $\obar\Bigl ( 0 ; \bigl (\max\{r_{1},r_{2}\} \bigr )^{2} \Bigr )$. This means that, for any of the
zeros $z_{j}$ of $p$, $|z_{j}|^{2} \leq \bigl ( \max\{r_{1},r_{2}\} \bigr )^{2}$, which concludes the proof of part (a) when $\ell$ is even.
When $\ell$ is odd, we multiply $p$ by $z$, which makes it a polynomial of even degree with an added zero at the origin and 
with $m=\ceil*{\frac{\ell}{2}}$. The coefficients $a_{j}$ of $zp(z)$ are
then as defined in the statement of the theorem, with $\azero = 0$. The latter is of no consequence and the proof of part (a) then follows from the even case. 
$ $ \newline
\newline {\bf (b)}
Here too, we first assume that $\ell$ is even. We then proceed similarly as in part (a) and let $\Delta_{x}$ be the block diagonal matrix 
with diagonal blocks ($x^{m}I,x^{m-1}I,\dots,xI$), where $I$ is the $2 \times 2$ identity matrix. We define
$\Phi_{x}(p) = \Delta_{x}^{-1} C^{2}_{S}(p) \Delta_{x}$, so that, with~(\ref{Phidef}), the Gershgorin column set of $\Phi_{x}(p)$ 
is the union of the three discs $O'_{1}(x) \equiv \obar(0;x)$,
$O'_{2}(x) \equiv \obar(\alpha;\sigma_{1}(x))$, and $O'_{3}(x) \equiv \obar(\beta;\sigma_{2}(x))$, where
$\sigma_{1}$ and $\sigma_{2}$ are defined by~(\ref{sigma1}) and~(\ref{sigma2}), respectively.
As $x$ increases, $O'_{1}(x)$ expands, while $O'_{2}(x)$ and $O'_{3}(x)$ contract. When $x = |\alpha| + \sigma_{1}(x)$, $O'_{1}(x)$ and $O'_{2}(x)$ are
tangent to one another and $O'_{2}(x) \subseteq O'_{1}(x)$. This occurs when $x=s_{1}$, where $s_{1}$ is the unique positive solution of
\bdis
x^{m} - |\alpha|x^{m-1} - \big ( |\vnmth | + |v_{n-4} | \bigr ) x^{m-2}  - \dots - \bigl ( |v_{3} | + |\vtwo | \bigr ) x 
- \dots - \lb |\vone | + |\vzero | \rb = 0 \; .
\edis
On the other hand, $O'_{1}(x)$ and $O'_{3}(x)$ are tangent to one another and $O'_{3}(x) \subseteq O'_{1}(x)$
when $x = |\beta| + \sigma_{2}(x)$, which happens when $x=s_{2}$, where $s_{2}$ is the unique positive solution of
\bdis
x^{m} - (|\beta|+|\gamma|) x^{m-1} - \big ( |\wnmth | + |w_{n-4} | \bigr ) x^{m-2}  - \dots - \bigl ( |w_{3} | + |\wtwo | \bigr ) x 
- \dots - \lb |\wone | + |\wzero | \rb = 0 \; .
\edis
From here on, the proof proceeds analogously to the proof of part (a) and we conclude that the Gershgorin set of $\Phi_{x}(p)$ is given by
$\obar\Bigl ( 0 ; \max\{s_{1},s_{2}\} \Bigr )$. Consequently, we obtain for the zeros $z_{j}$ of $p$ that
$|z_{j}|^{2} \leq \max\{s_{1},s_{2}\}$, which concludes the proof of part (b) for even $\ell$. When $\ell$ is odd, we consider $zp(z)$ instead of $p(z)$
and the proof follows from the even case, analogously as in part (a). \qquad \qed

One of the advantages of the triangularization of $\Ammo$ when applying Gershgorin's theorem is made clear by the proof of part (a), where $O_{2}(x)$ would 
otherwise not necessarily contract with increasing $x$ if the (2,1)-element of $\Ammo$ were nonzero, as it would be multiplied by $x$ after the similarity 
transformation. 

Although the solution of the real equations in this theorem requires only a fraction of the computational effort required to compute all the zeros of $p$,
it is worth mentioning that this can be carried out very efficiently. Once, e.g., $r_{1}$ is computed, the sign of $\psi_{2}(r_{1})$ immediately determines
if it is larger than $r_{2}$ or not. If it is, $\psi_{2}(x)=0$ need not be solved. An analogous situation exists for part (b). Moreover, it is computationally
less expensive to solve two equations of degree $n/2$ than just one of degree $n$.

Although it would lead us too far, more polynomial inclusion regions along the lines of the ones obtained in~\cite{MelmanTwin} or \cite[Corollary 8.2.3]{RS} 
can be derived here as well. Such regions are obtained when one of the discs centered at $\alpha$ or $\beta$ is allowed to absorb the one centered at the origin. 
In addition, other special values for the parameter $x$ might also be considered, such as, e.g., values for which two discs have the same radius.  

%
%
                                                            \section{Pellet-like results}
                                                             \label{PelletResults}                                            
%
%

It is sometimes possible to isolate one zero of a polynomial by applying Gershgorin's theorem to $C(p)$ (see, e.g. \cite{MelmanTwin}), and 
a similar approach can be applied here, although here it can lead to the isolation of one or two squares of zeros of the polynomial (and therefore to the isolation
of the zeros themselves). This is reminiscent of special cases of Pellet's theorem, whence the title of this section, although it does provide a smaller       
inclusion region than those provided by Pellet's theorem, which lead to discs, annuli, and (infinite) complements of discs. 
Very often, however, it is the very ability to isolate zeros that is 
important, rather than the size of the inclusion region. We formulate two theorems, based on the similarity transformations of $C^{2}(p)$ introduced in 
Section~\ref{Preliminaries}. To make them self-contained, their statements include some previously defined quantities.
%
%
\begin{theorem}
\label{TGP}                       
For a polynomial $p(z)= z^{\ell} + b_{\ell-1} z^{\ell-1} + \dots + b_{1} z + b_{0}$ with complex coefficients, $\ell \geq 3$, 
define $m = \ceil*{\frac{\ell}{2}}$, $n=2m$, and the complex numbers $a_{j}$ as follows:
\begin{eqnarray*}
& & a_{j} = b_{j} \;\; \text{($\ell$ is even and $j=0,\dots,n-1$),}  \\
& & a_{0} = 0 \; \text{and} \; a_{j} = b_{j-1} \;\; \text{($\ell$ is odd and $j=1,\dots,n-1$).}  
\end{eqnarray*}
Furthermore, define 
\bdis
A_{0} = \begin{pmatrix} -\azero & \anmo \azero \\ -\aone & \anmo \aone -\azero \end{pmatrix} \; , 
\;\; \text{and} \;\;
A_{j} = \begin{pmatrix} -a_{2j} & \anmo a_{2j} -a_{2j-1} \\ -a_{2j+1} & \anmo\ajpo -a_{2j} \\ \end{pmatrix} 
\;\; \text{for $j=1,...,m-1$,}
\edis
and let $S$ be a nonsingular matrix such that 
$            
S^{-1}\Ammo S = \begin{pmatrix} \alpha & \gamma \\ 0 & \beta \\ \end{pmatrix} \; ,
$           
where $\alpha, \beta, \gamma \in \complex$. Define the vectors $v,w \in \complex^{n-2}$ as                      
\bdis
\begin{pmatrix} v_{2j} & w_{2j} \\ v_{2j+1} & w_{2j+1} \\ \end{pmatrix} = S^{-1} \Aj S  
\;\; \text{for $j=0,...,m-2$} \; .
\edis
For $x > 0$, let
\begin{eqnarray*}
\rho_{1}(x)   & = & \dfrac{|\vnmth |}{x} + \dfrac{|v_{n-4} |}{x^{2}} + \dots + \dfrac{|\vone |}{x^{n-3}} + \dfrac{|\vzero |}{x^{n-2}} \; ,    \\ 
\rho_{2}(x)   & = & \dfrac{|\gamma|}{x} + \dfrac{|\wnmth |}{x^{2}} + \dfrac{|w_{n-4} |}{x^{3}} + \dots + \dfrac{|\wone |}{x^{n-2}} + \dfrac{|\wzero |}{x^{n-1}}  \; , \\ 
\chi_{1}(x)   & = & x^{n} - |\alpha|x^{n-2} + |\vnmth |x^{n-3} + \dots + |\vone | x + |\vzero | \; ,        \\    
\chi_{2}(x)   & = & x^{n+1} - |\beta|x^{n-1} + |\gamma| x^{n-2} + |\wnmth |x^{n-3} + \dots + |\wone | x + |\wzero |   \; ,     
\end{eqnarray*}
and
\begin{eqnarray*}
\sigma_{1}(x) & = & \dfrac{|\vnmth | + |v_{n-4} |}{x} + \dfrac{|v_{n-5}| + |v_{n-6}|}{x^{2}} + \dots + \dfrac{|\vone | + |\vzero |}{x^{m-1}} \; ,  \\ 
\sigma_{2}(x) & = & |\gamma| + \dfrac{|\wnmth | + |w_{n-4} |}{x} + \dfrac{|w_{n-5}| + |w_{n-6}|}{x^{2}} + \dots + \dfrac{|\wone | + |\wzero |}{x^{m-1}} \; ,  \\
\omega_{1}(x) & = & x^{m} - |\alpha|x^{m-1} + \big ( |\vnmth | + |v_{n-4} | \bigr ) x^{m-2}  
+ \dots + \bigl ( |v_{3}| + |v_{2}| \bigr ) x + \dots + \lb |\vone | + |\vzero | \rb  \; ,   \\
\omega_{2}(x) & = & x^{m} - (|\beta|-|\gamma|) x^{m-1} + \big ( |\wnmth | + |w_{n-4} | \bigr ) x^{m-2}  
+ \dots + \bigl ( |w_{3}| + |w_{2}| \bigr ) x + \dots + \lb |\wone | + |\wzero | \rb  \; .
\end{eqnarray*}
Define for $x>0$: 
\begin{eqnarray*}
& & O_{1}(x) \equiv \obar(0;x^{2}) \; , \; O_{2}(x) \equiv \obar(\alpha;\rho_{1}(x)) \; , \; O_{3}(x) \equiv \obar(\beta;\rho_{2}(x)) \; ,  \\
& & O'_{1}(x) \equiv \obar(0;x) \; , \; O'_{2}(x) \equiv \obar(\alpha;\sigma_{1}(x)) \; ; \; O'_{3}(x) \equiv \obar(\beta;\sigma_{2}(x)) \; .
\end{eqnarray*}
If $\chi_{1}(x)=0$ has positive solutions $t_{1}$ and $t_{2}$, set $I_{1}=[t_{1},t_{2}]$, otherwise $I_{1}=\emptyset$.
If $\chi_{2}(x)=0$ has positive solutions $u_{1}$ and $u_{2}$, set $I_{2}=[u_{1},u_{2}]$, otherwise $I_{2}=\emptyset$.
If $\omega_{1}(x)=0$ has positive solutions $\tildet_{1}$ and $\tildet_{2}$, set $J_{1}=[\tildet_{1},\tildet_{2}]$, otherwise $J_{1}=\emptyset$.
If $\omega_{2}(x)=0$ has positive solutions $\tildeu_{1}$ and $\tildeu_{2}$, set $J_{2}=[\tildeu_{1},\tildeu_{2}]$, otherwise $J_{2}=\emptyset$.
Then the following holds.                    
\newline {\bf (a1)} If $I_{1} \cap I_{2} \neq \emptyset$, set 
$\mu_{1} = \max\{t_{1},u_{1}\}$ and $\mu_{2} = \min\{t_{2},u_{2}\}$. Then $\ell-2$ of $p$'s zeros lie in $\obar(0;\mu_{1})$,
while the union $O_{2}(\rho_{1}(\mu_{2})) \cup O_{3}(\rho_{2}(\mu_{2}))$ contains the squares of the two remaining zeros, i.e., $p$ has no zeros 
with a modulus between $\mu_{1}$ and $\left ( \min \bigl \{ |\alpha| - \rho_{1}(\mu_{2}), |\beta| - \rho_{2}(\mu_{2}) \bigr \} \right )^{1/2}$.
If, in addition, $|\alpha - \beta | > \rho_{1}(\mu_{2}) + \rho_{2}(\mu_{2})$, then $O_{2}(\mu_{2})$ and $O_{3}(\mu_{2})$ 
each contain one square 
of a zero of $p$.
\newline {\bf (a2)} If $I_{1} \cap I_{2} = \emptyset$, we have the following.
\begin{itemize}
\item 
If $I_{1} \neq \emptyset$ and $|\alpha-\beta| > \rho_{1}(t_{2}) + \rho_{2}(t_{2})$, then 
the squares of $\ell-1$ of $p$'s zeros 
are contained in the union $O_{1}(t_{2}) \cup O_{3}(t_{2})$, while the remaining square of a zero lies in $O_{2}(t_{2})$.
\item 
If $I_{2} \neq \emptyset$ and $|\alpha-\beta| > \rho_{1}(u_{2}) + \rho_{2}(u_{2})$, then 
the squares of $\ell-1$ of $p$'s zeros 
are contained in the union $O_{1}(u_{2}) \cup O_{2}(u_{2})$, while the remaining square of a zero lies in $O_{3}(u_{2})$.
\end{itemize}
{\bf (b1)} If $J_{1} \cap J_{2} \neq \emptyset$, set 
$\tildemu_{1} = \max\{\tildet_{1},\tildeu_{1}\}$ and $\tildemu_{2} = \min\{\tildet_{2},\tildeu_{2}\}$. Then $\ell-2$ of $p$'s zeros lie in 
$\obar(0;\sqrt{\tildemu_{1}})$,
while the union $O'_{2}(\sigma_{1}(\tildemu_{2})) \cup O'_{3}(\sigma_{2}(\tildemu_{2}))$ contains the squares of the two remaining zeros, i.e., $p$ has no zeros 
with a modulus between $\sqrt{\tildemu_{1}}$ and $\left ( \min \bigl \{ |\alpha| - \sigma_{1}(\tildemu_{2}), |\beta| - \sigma_{2}(\tildemu_{2}) \bigr \} \right )^{1/2}$.
If, in addition, $|\alpha - \beta | > \sigma_{1}(\tildemu_{2}) + \sigma_{2}(\tildemu_{2})$, then $O'_{1}(\tildemu_{2})$ and $O'_{2}(\tildemu_{2})$ 
each contains one square of a zero of $p$.
\newline {\bf (b2)} If $I_{1} \cap I_{2} = \emptyset$, we have the following.
\begin{itemize}
\item 
If $J_{1} \neq \emptyset$ and $|\alpha-\beta| > \sigma_{1}(\tildet_{2}) + \sigma_{2}(\tildet_{2})$, then 
the squares of $\ell-1$ of $p$'s zeros 
are contained in the union $O'_{1}(\tildet_{2}) \cup O'_{3}(\tildet_{2})$, while the remaining square of a zero lies in $O'_{2}(\tildet_{2})$.
\item 
If $J_{2} \neq \emptyset$ and $|\alpha-\beta| > \sigma_{1}(\tildeu_{2}) + \sigma_{2}(\tildeu_{2})$, then 
the squares of $\ell-1$ of $p$'s zeros 
are contained in the union $O'_{1}(\tildeu_{2}) \cup O'_{2}(\tildeu_{2})$, while the remaining square of a zero lies in $O'_{3}(\tildeu_{2})$.
\end{itemize}
\end{theorem}
\noindent \prf
\newline {\bf (a1)} When $\ell$ is even, $m=\ell/2$, $n=\ell$, and $a_{j}=b_{j}$.
In the Gershgorin column set for $F_{x}(p)$, defined by~(\ref{Fdef}), $O_{1}(x)$ is disjoint from the other two discs if 
\bdis
|\alpha| > x^{2} + \rho_{1}(x) \Longleftrightarrow \chi_{1} < 0  \;\; \text{and} \;\; |\beta| > x^{2} + \rho_{2}(x) \Longleftrightarrow \chi_{2} < 0 \; .
\edis
This can only happen if $I_{1} \cap I_{2}$ is not empty, in which case $I_{1} \cap I_{2} = [\mu_{1},\mu_{2}]$, where $\mu_{1}$ and $\mu_{2}$
are defined in the statement of the theorem. By Gershgorin's theorem, $O_{1}(x)$ then contains $n-2$ squares of zeros of $p$, while 
$O_{2}(x) \cup O_{3}(x)$ contains the remaining two for any $x$ satisfying $\mu_{1} < x < \mu_{2}$. This is therefore true for the intersection
of all these Gershgorin sets as $x$ runs from $\mu_{1}$ to $\mu_{2}$, which is given by the disjoint sets $O_{1}(\mu_{1}^{2})$ and 
$O_{2}(\mu_{2}) \cup O_{3}(\mu_{2})$. As a consequence, and because $O_{2}(x)$ and $O_{3}(x)$ are centered at $\alpha$ and $\beta$, respectively,
no square of the zeros of $p$ can have a modulus between $\mu_{1}^{2}$ and 
$\min \{ |\alpha| - \rho_{1}(\mu_{2}), |\beta| - \rho_{2}(\mu_{2}) \}$. 
If $|\alpha - \beta | > \rho_{1}(\mu_{2}) + \rho_{2}(\mu_{2})$, then $O_{2}(\mu_{2})$ and $O_{3}(\mu_{2})$ are also disjoint from each other, and by
Gershgorin's theorem must each contain a square of a zero of $p$. For odd $n$, we consider $zp(z)$ instead of $p(z)$, as in the proof of Theorem~\ref{TGC},
and then proceed as in the even case. This proves part (a1).
\newline {\bf (a2)} Assume that $\ell$ is even.
If $I_{1} \cap I_{2} = \emptyset$, then there are no values of $x$ for which both $O_{2}(x)$ and $O_{3}(x)$ are disjoint from $O_{1}(x)$.
When $I_{1} \neq \emptyset$, then the smallest radius $O_{2}(x)$ can have while being disjoint from $O_{1}$ is $x=\rho_{2}(t_{2})$. If $O_{2}(t_{2})$
and $O_{3}(t_{2})$ are disjoint, i.e., when $|\alpha-\beta| > \rho_{1}(t_{2}) + \rho_{2}(t_{2})$, then by Gershgorin's theorem, 
$O_{2}(t_{2})$ contains exactly one square of a zero of $p$, while the other $n-1$ are contained in $O_{1}(t_{2}) \cup O_{3}(t_{2})$. 
The proof of the analogous situation, where the roles of $I_{1}$ and $I_{2}$ are switched, then also follows. The case where $\ell$ is odd is treated as before.
This proves part (a2).
\newline {\bf (b1)} and {\bf (b2)} The proof of parts (b1) and (b2) follows the exact same pattern as that of parts (a1) and (a2), except that here
the radius of $O'_{1}(x)$ is $x$ instead of $x^{2}$. All other aspects are analogous. This concludes the proof of the theorem. \qed

We note that $|\beta| > |\gamma|$ is a necessary condition for $\omega_{2}(x)=0$ to have positive solutions.

\begin{figure}[H]
\begin{center}
\raisebox{0ex}{\includegraphics[width=0.45\linewidth]{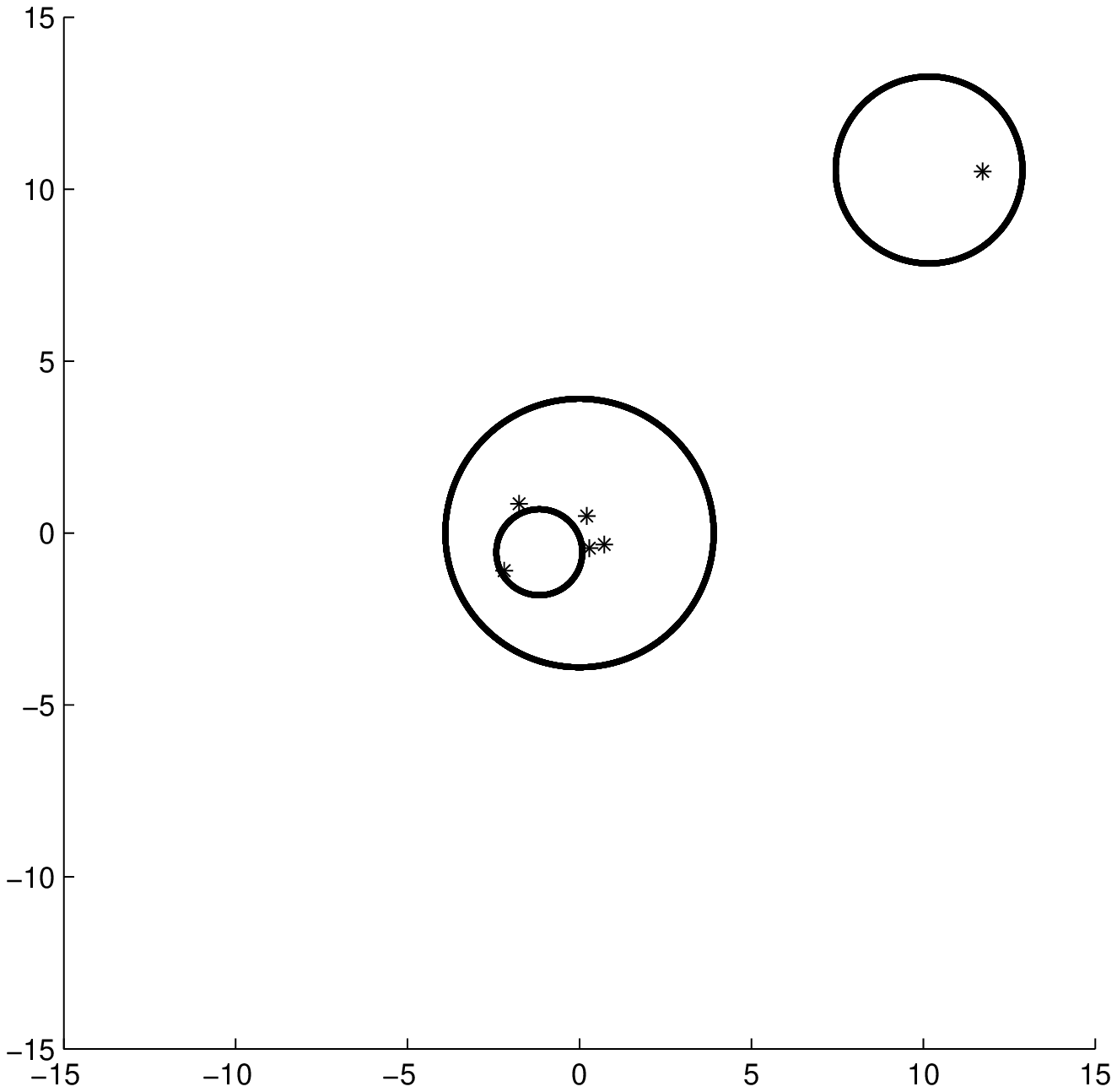}}
\qquad  
\raisebox{0ex}{\includegraphics[width=0.45\linewidth]{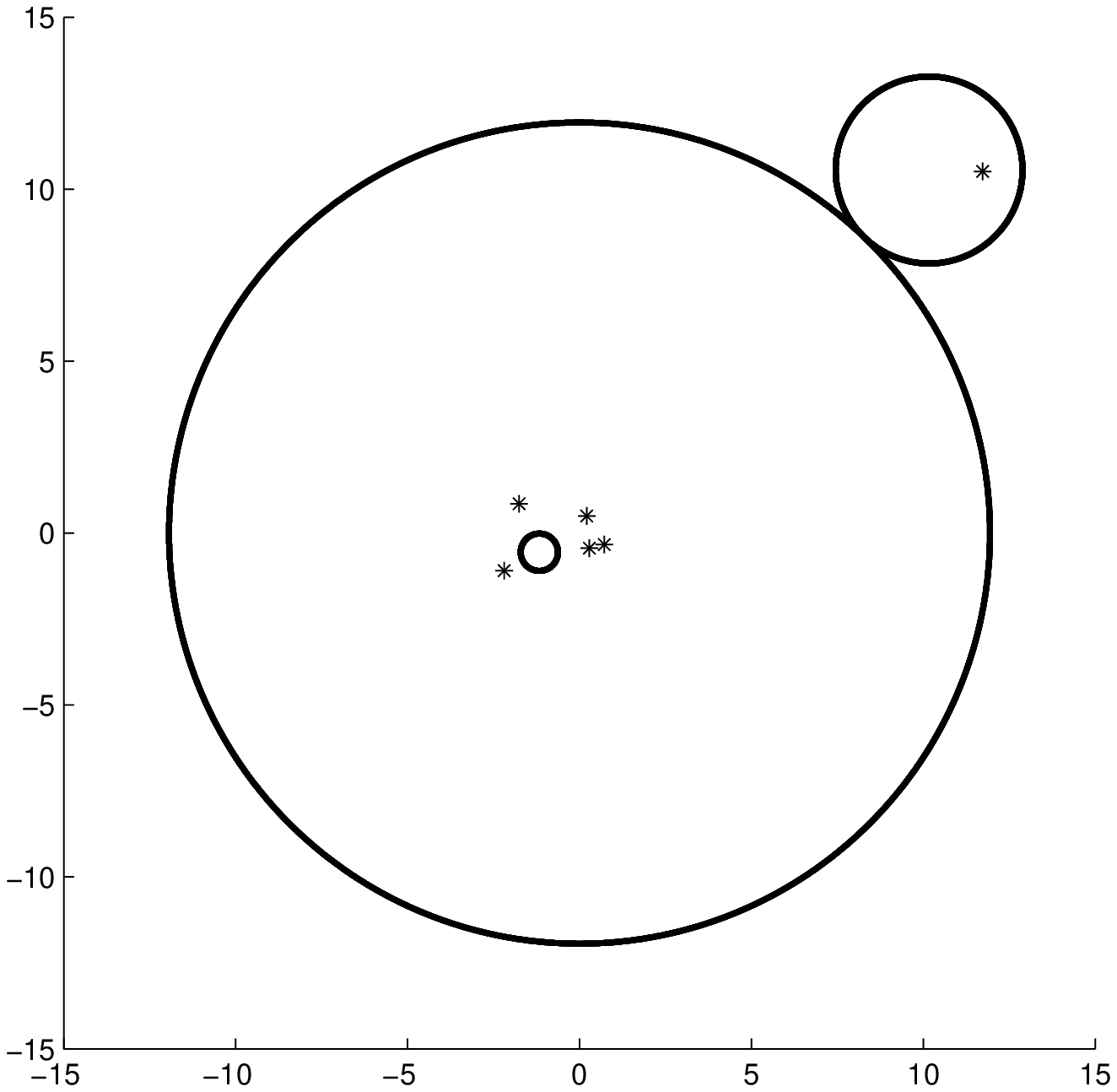}}
\caption{Enhancement of case (a2) with $u^{*} < u_{1}$.}                                          
\label{case_a2_detached}
\end{center}
\end{figure}

\begin{figure}[H]
\begin{center}
\raisebox{0ex}{\includegraphics[width=0.45\linewidth]{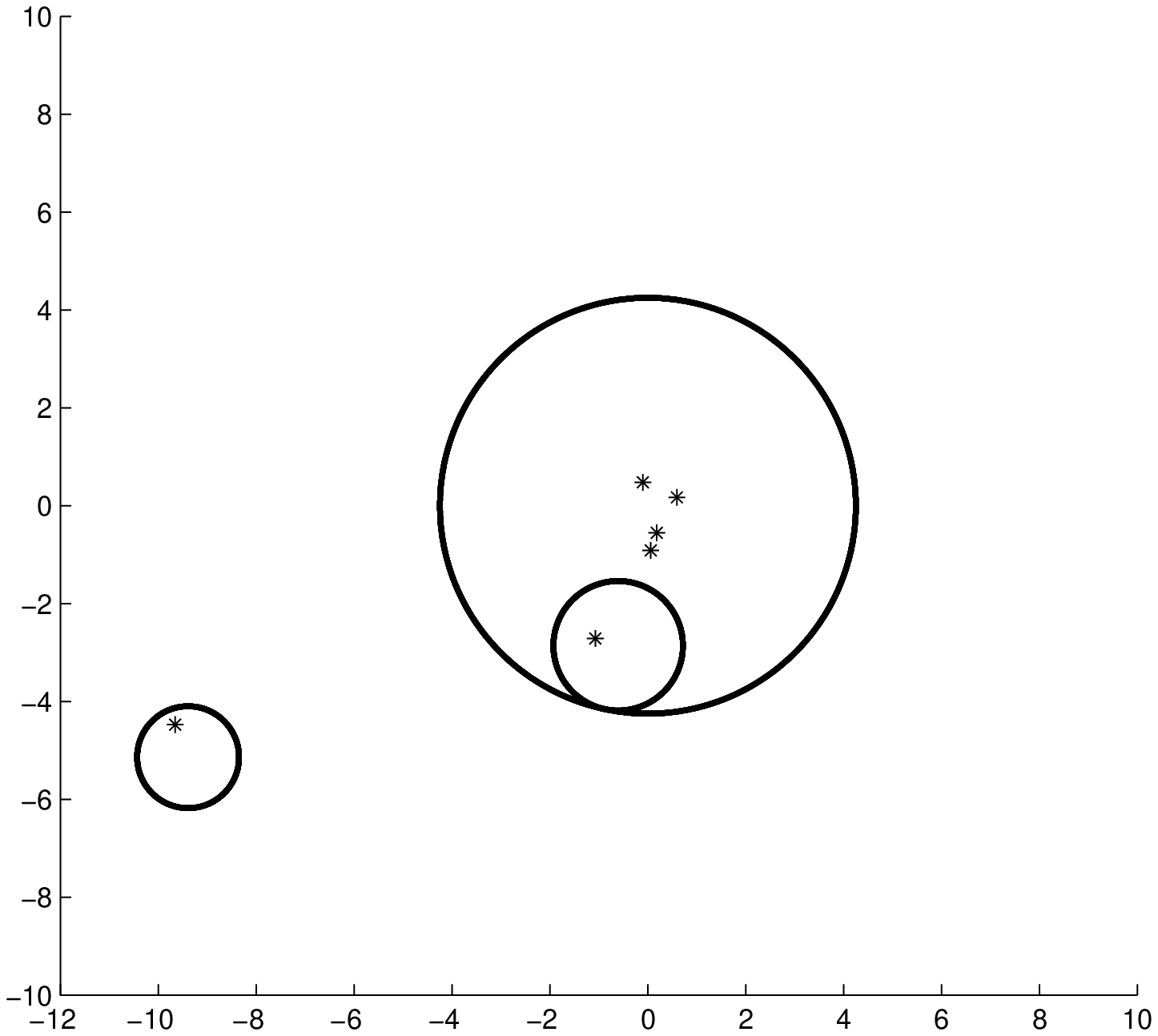}}
\qquad  
\raisebox{0ex}{\includegraphics[width=0.45\linewidth]{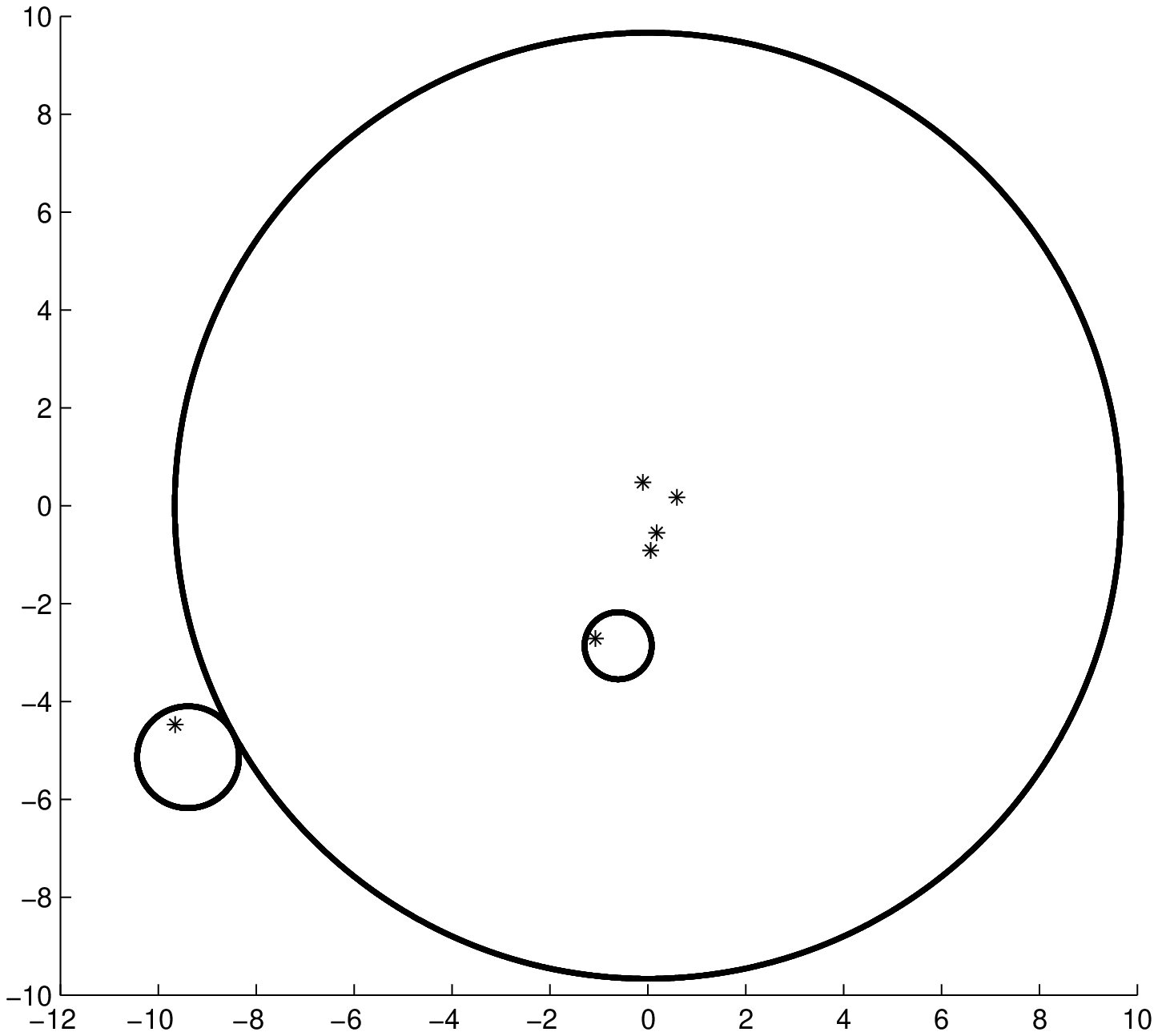}}
\caption{Enhancement of case (a2) with $u^{*} > u_{1}$.}                                          
\label{case_a2_tangent}  
\end{center}
\end{figure}

The previous theorem, which focuses on the isolation of zeros, but does not try to optimize the inclusion sets, can be enhanced in several ways, 
depending on the situation. As an example, consider the case (a2), where $I_{2} \neq \emptyset$ and
$|\alpha - \beta| > \rho_{1}(u_{2}) + \rho_{2}(u_{2})$. If $O_{2}(u_{2})$ happens to be contained in $O_{1}(u_{2})$, then there exists a value 
$u^{*} < u_{2}$ for which $O_{2}(u^{*})$ still lies inside $O_{1}(u^{*})$, but is tangent to it. In that case, the Gershgorin column set
of $F_{x}(p)$, given by
the union of the two disjoint sets $O_{3}(u_{2})$ and $O_{1}\left (\max\{ u_{1}, u^{*} \} \right )$, is smaller than the one in the theorem.
In fact, in this case, even when $O_{2}(u_{2})$ does not lie inside $O_{1}(u_{2})$, $O_{1}(u) \cup O_{2}(u)$ for some values $u < u_{2}$ may be 
smaller than the corresponding set in the theorem. Other cases may be similarly improved.
Figure~\ref{case_a2_detached} shows an example of case (a2) as we just described with $u^{*} < u_{1}$: on the left is the enhanced Gershgorin set, while
the Gershgorin set from the theorem is shown on the right. It was obtained for the polynomial $z^{6}+(4+3i)z^{5}+(3-i)z^{4}r+4+2i)z^{3}-(3+2i)z^{2}-3z+4+i$. 
Figure~\ref{case_a2_tangent} shows a similar example, where this time $u^{*} > u_{1}$, obtained 
for the polynomial $z^{6}-2iz^{5}+(3+4i)z^{4}+(3+i)z^{3}-(2+i)z^{2}+2iz+2+i$. The asterisks in the figures indicate the squares of the zeros of the polynomials. 

%
%
%
%

%
%
%
%

In the introduction, we mentioned that Theorem~\ref{CauchyTheorem} can be obtained by applying Gershgorin's theorem to a similarity transformation of $C(p)$, 
and it can be shown that similarly applying the Block Gershgorin theorem to $C^{2}_{S}(p)$ is equivalent to applying its matrix version, namely, 
Theorem~\ref{MatrixCauchyTheorem}, to the matrix polynomial $P_{S}$, defined by~(\ref{PSdef}), when $\Ammo$ is diagonalizable.
However, for Pellet's theorem, applying the Block Gershgorin theorem to $C^{2}_{S}(p)$ leads to a more subtle result than its matrix version,
which is the next theorem. For this theorem we will assume $\Ammo$ to be diagonalizable.

%
%
\begin{theorem}
\label{TMGP}                       
For a polynomial $p(z)= z^{\ell} + b_{\ell-1} z^{\ell-1} + \dots + b_{1} z + b_{0}$ with complex coefficients, $\ell \geq 3$, 
define $m = \ceil*{\frac{\ell}{2}}$, $n=2m$, and the complex numbers $a_{j}$ as follows:
\begin{eqnarray*}
& & a_{j} = b_{j} \;\; \text{($\ell$ is even and $j=0,\dots,n-1$),}  \\
& & a_{0} = 0 \; \text{and} \; a_{j} = b_{j-1} \;\; \text{($\ell$ is odd and $j=1,\dots,n-1$).}  
\end{eqnarray*}
Furthermore, define 
\bdis
A_{0} = \begin{pmatrix} -\azero & \anmo \azero \\ -\aone & \anmo \aone -\azero \end{pmatrix} \; , 
\;\; \text{and} \;\;
A_{j} = \begin{pmatrix} -a_{2j} & \anmo a_{2j} -a_{2j-1} \\ -a_{2j+1} & \anmo\ajpo -a_{2j} \\ \end{pmatrix} 
\;\; \text{for $j=1,...,m-1$,}
\edis
and let $\Ammo$ be diagonalizable by a matrix $S$ such that 
$            
S^{-1}\Ammo S = \begin{pmatrix} \alpha & 0 \\ 0 & \beta \\ \end{pmatrix} \; ,
$           
where $\alpha, \beta \in \complex$ and $|\alpha| \leq |\beta|$. 
Define the $2 \times 2$ matrices $T_{j}=S^{-1} \Aj S$ ($j=0,\dots,m-1$) and the complex vectors $v$ and $w$
as
\bdis
\begin{pmatrix} v_{2j} & w_{2j} \\ v_{2j+1} & w_{2j+1} \\ \end{pmatrix} = T_{j}  \;\; \text{for $j=0,...,m-2$} \; .
\edis
For $x > 0$, let
\begin{eqnarray*}
\tau(x) & = & \dfrac{||T_{m-2}||}{x} + \dfrac{||T_{m-3}||}{x^{2}} + \dots + \dfrac{||T_{0}||}{x^{m-1}} \; ,  \\
\Omega_{1}(x) & = & x^{m} - |\alpha|x^{m-1} + ||T_{m-2}|| x^{m-2}  + \dots + ||T_{1}|| x + \dots + ||T_{0}||  \; ,   \\
\Omega_{2}(x) & = & x^{m} - |\beta|x^{m-1} + ||T_{m-2}|| x^{m-2}  + \dots + ||T_{1}|| x + \dots + ||T_{0}||  \; ,                
\end{eqnarray*}
and define 
\bdis
O''_{1}(x) \equiv \obar(0;x) \; , \; O''_{2}(x) \equiv \obar(\alpha;\tau(x)) \; ; \; O''_{3}(x) \equiv \obar(\beta;\tau(x)) \; .
\edis
If $\Omega_{1}(x)=0$ has positive solutions $x_{1}$ and $x_{2}$, set $K_{1}=[x_{1},x_{2}]$, otherwise $K_{1}=\emptyset$.
If $\Omega_{2}(x)=0$ has positive solutions $y_{1}$ and $y_{2}$, set $K_{2}=[y_{1},y_{2}]$, otherwise $K_{2}=\emptyset$.
Then the following holds.                    
\newline {\bf (a)} If $K_{1} \neq \emptyset$, then $K_{2} \neq \emptyset$, $K_{1} \subseteq K_{2}$, and 
$\ell-2$ of $p$'s zeros lie in $O''_{1}(\sqrt{x_{1}})$,
while the union $O_{2}(\tau(x_{2})) \cup O_{3}(\tau(x_{2}))$ contains the squares of the two remaining zeros, i.e., $p$ has no zeros 
with a modulus between $\sqrt{x_{1}}$ and $\bigl ( |\alpha| - \tau(x_{2}) \bigr )^{1/2}$.
If, in addition, $|\alpha - \beta | > 2 \tau(x_{2})$, then $O''_{2}(x_{2})$ and $O''_{3}(x_{2})$ 
each contains one square 
of a zero of $p$.
\newline {\bf (b)} If $K_{1} = \emptyset$ and $K_{2} \neq \emptyset$, then if 
$|\alpha-\beta| > 2\tau(y_{2})$, 
the squares of $\ell-1$ of $p$'s zeros 
are contained in the union $O''_{1}(y_{2}) \cup O''_{2}(y_{2})$, while the remaining square of a zero lies in $O''_{3}(y_{2})$.
\end{theorem}
\prf The proof is similar to the proof of Theorem~\ref{TGP} with minor differences. Assume first that $\ell$ is even. 
Here we apply Theorem~\ref{BlockGershgorinTheorem}, the Block
Gershgorin theorem, to $\Phi_{x}(p)$, defined in~(\ref{Phidef}). As we showed in Section~\ref{Preliminaries}, this produces the block Gershgorin 
column set $G_{1} \cup G_{2}$, with the sets $G_{j}$ as defined in the statement of Therorem~\ref{BlockGershgorinTheorem}. There we saw that
$G_{1} = O''_{1}(x)$ and that $G_{2} = O''_{2}(x) \cup O''_{3}(x)$, the union of two discs centered at $\alpha$ and $\beta$, respectively, with identical 
radii $\tau(x)$.
$O''_{1}(x)$ is disjoint from the other two discs, if 
\bdis
|\alpha| > x + \tau(x) \Longleftrightarrow \Omega_{1}(x) < 0  \;\; \text{and} \;\; |\beta| > x + \tau(x) \Longleftrightarrow \Omega_{2}(x) < 0 \; .
\edis
Clearly, because $|\alpha| \leq |\beta|$, if $\Omega_{1}(x)=0$ has two positive solutions $x_{1}$ and $x_{2}$, then $\Omega_{2}(x)=0$
also has two positive solutions $y_{1}$ and $y_{2}$, with $y_{1} \leq x_{1} \leq x_{2} \leq y_{2}$. From here on, the proof follows analogously to that of 
Theorem~\ref{TGP}. The case when $\ell$ is odd is treated analogously. \qed

Similar enhancements of this theorem can be obtained as for Theorem~\ref{TGP}.

\begin{figure}[H]
\begin{center}
\raisebox{0ex}{\includegraphics[width=0.5\linewidth]{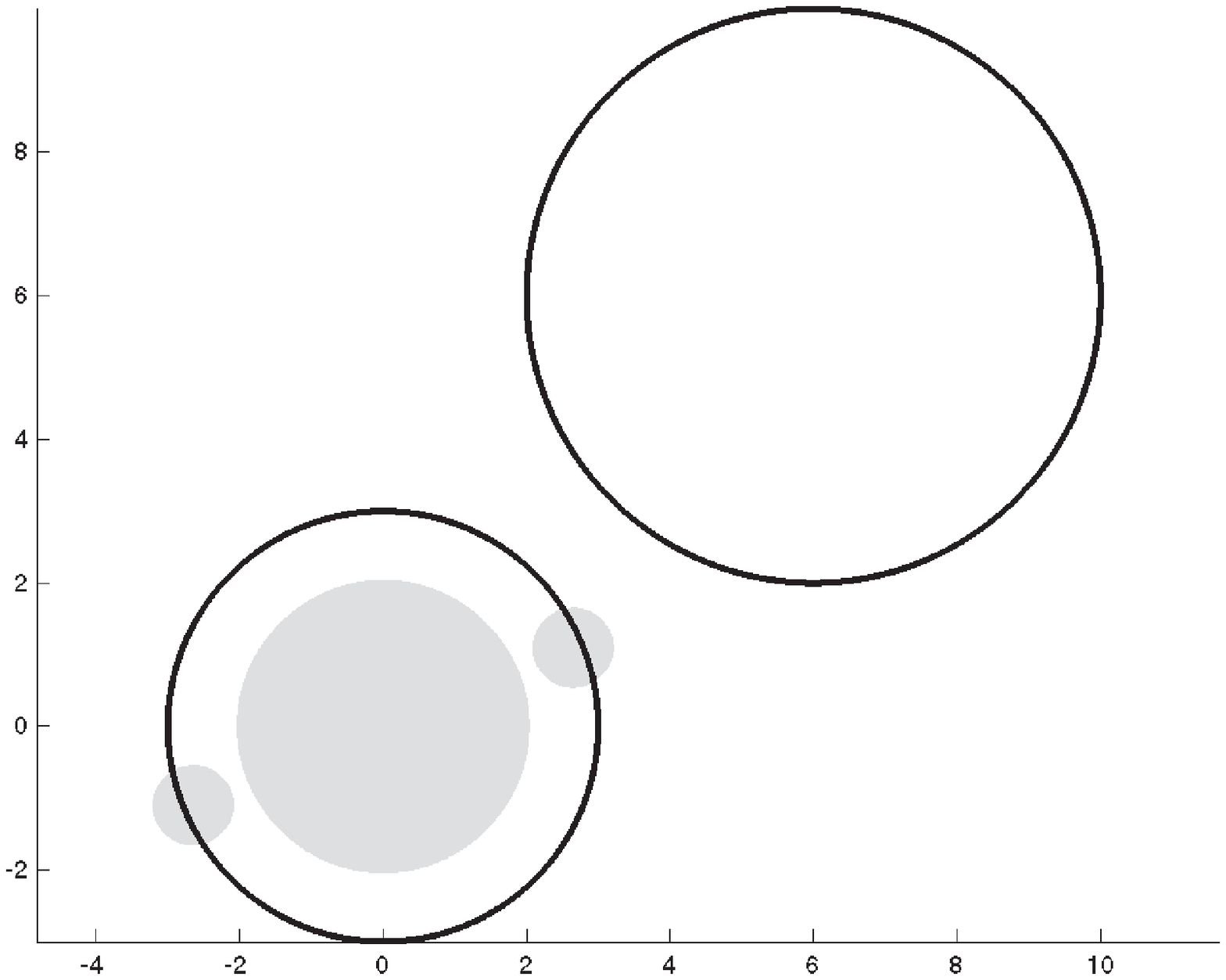}}
\caption{Inclusion regions for the squares of the zeros (circles) and the zeros themselves (shaded).}                
\label{Cassini}                   
\end{center}
\end{figure}

The theorems in this section derive inclusion regions that are sometimes given in terms of the squares of the zeros of $p$. These results are easily
translated to bounds on and gaps between the moduli of the zeros of $p$. However, the inclusion sets themselves are slightly more complicated.
If a disc, centered at $c$ with radius $R$, contains the squares of the zeros of a polynomial, then any zero $z$ satisfies $|z^{2} - c| \leq R$,
which implies that $|z+\sqrt{c}||z-\sqrt{c}| \leq R$, i.e., they lie in a region bounded by an oval of Cassini with foci $\pm \sqrt{c}$.
As illustration let us consider a situation where the squares of the zeros of a 
polynomial are contained in the union of a disc, centered at the origin with radius $R_{1}$, and another disc, centered at $c$ with radius $R_{2}$. 
Figure~\ref{Cassini} shows such discs, containing the squares of the zeros with $c=6+6i$, $R_{1}=4$, and $R_{2}=3$, while 
the corresponding inclusion region for the zeros themselves - the union of a disc and an oval of Cassini (consisting of two loops
because the disc centered at $c$ is bounded away from the origin) - are shaded in gray.

We conclude by pointing out that a similar approach has the potential to improve analogous results for matrix polynomials, especially when the matrix 
coefficients are of moderate size compared to the degree of the polynomial, in which case one can argue that the corresponding block companion 
matrices also have a diagonal mostly composed of zeros; squaring them may lead to smaller inclusion regions there as well.
%
%
%
%

%
%
                                                   \section{Numerical comparisons.} 
                                                   \label{NumericalComparisons} 
%
%
In this section we illustrate our results numerically, while comparing them to the classical Theorems~\ref{CauchyTheorem} and \ref{PelletTheorem}.
To do so, we created sets of random polynomials, to which, for the Cauchy-like results of Section~\ref{CauchyResults}, we applied
Theorem~\ref{CauchyTheorem}, Theorem~\ref{MatrixCauchyTheorem} applied to the matrix polynomial $P_{S}$ defined by~(\ref{PSdef}),
and parts (a) and (b) of Theorem~\ref{TGC}, while for the  
Pellet-like results of Section~\ref{PelletResults} we compared Theorem~\ref{PelletTheorem} with $k=1,2$, Theorem~\ref{MatrixPelletTheorem} 
applied to the aforementioned matrix polynomial $P_{S}$ with $k=1$, all parts of Theorem~\ref{TGP}, and Theorem~\ref{TMGP}. 
For the Cauchy-like results, we compared the averages of the ratio of the upper bounds and
the modulus of the largest zero, i.e., the closer this number is to one, the better the bound, and we also recorded the number of times each method gave the
best upper bound on the moduli of the zeros. 
For the Pellet-like results, we compared the number of times zeros (or squares of zeros) could be isolated from the others for each method,
which is generally the most important use of these methods.
We chose the Euclidean norm (2-norm) when applying Theorem~\ref{BlockGershgorinTheorem} and the diagonalizing matrix $S$ was chosen to have normalized 
columns. No polynomials were generated where $\Ammo$ was not diagonalizable, and its eigenvalues $\alpha$ and $\beta$ as they appear in all our results
were ordered so that $|\alpha| \leq |\beta|$. No special choice of coefficients was made in Sets~1-2, but Sets~3-4 exhibit 
specific relative magnitudes of some coefficients to better illustrate the advantages of our methods.  

The sets of polynomials are defined below, with $n$ indicating the degree of the polynomials.
\newline {\bf Set 1:} n=20, the leading coefficient is one and the other coefficients have real and imaginary parts that are uniformly randomly distributed
on the interval $[-2,2]$. 
\newline {\bf Set 2:} n=20, the leading coefficient is one and the other coefficients have real and imaginary parts that are uniformly randomly distributed
on the interval $[-4,4]$. 
\newline {\bf Set 3:} n=20, the first four coefficients are 1,2,6,2 and the other coefficients have real and imaginary parts that are uniformly randomly distributed
on the interval $[-4,4]$. 
\newline {\bf Set 4:} n=20, the first four coefficients are 1,2,8,2 and the other coefficients have real and imaginary parts that are uniformly randomly distributed
on the interval $[-4,4]$. 

For each set we generated 1000 random polynomials and collected the results in Table~\ref{table1} for the Cauchy-like methods and in 
Table~\ref{table2} for the Pellet-like methods. In Table~\ref{table1}, the methods are listed across the top, and in an entry of the form $\gamma / j$,
$\gamma$ is the average ratio of the upper bound to the modulus of the largest zero, while $j$ is the number of times that a particular method
delivered the best ratio. It is clear from these results that the classical result by Cauchy (Theorem~\ref{CauchyTheorem}) is almost always worse than the 
other methods.
In Table~\ref{table2}, the methods are listed as in Table~\ref{table1}, and in an entry of the form $i / j / k$, $i$ is the number of times two zeros 
could be isolated not only from the $n-2$ remaining ones, but also from each other, $j$ is the number of times two zeros could be isolated from the other 
$n-2$ ones, but not from each other, and $k$ is the number of times a single zero could be isolated. For Pellet's theorem (first column), 
in an entry of the form $i / j$, $i$ and $j$ are the number of times one and two zeros could be isolated, respectively, from the remaining zeros. 
For Pellet's theorem's matrix version (second column), we listed the number of times two zeros could be isolated from the remaining $n-2$.

Our Pellet-like methods appear to be more sensitive, i.e., better able to isolate zeros, and do not seem to require as large a difference between the 
magnitudes of appropriate coefficients as is the case for Pellet's theorem. Moreover, the zero inclusion regions defined by Pellet's theorem are cruder 
than the results derived here.
The difference with Pellet's theorem can be quite dramatic, as for Set~1, where Pellet's Theorem was able to isolate zeros in only one case,
compared to more than 140 cases for our methods, and also Set~2, where it was able to isolate zeros for only 119 cases as opposed to more than 580 for our methods. 
A similar observation holds in the case of the isolation of two zeros for Set~3 and Set~4. 

We observed that the results for all sets of polynomials did not seem sensitive to the degree of the polynomial, delivering very similar results when the 
degree was doubled nor is there any appreciable difference between even and odd degrees. 
However, they are sensitive to the range of the real and complex parts of the generated random polynomials for the Pellet-like
results: the larger the range, the better the results, as it caused larger differences between the magnitudes of the coefficients, thereby increasing the 
likelihood that zeros can be separated. The effect of this on the Cauchy-like results was not significant.

We note that in~\cite{Melman_MatPol} the matrix version of Pellet's theorem was applied to $C^{2}(p)$, but not $C_{S}^{2}(p)$, 
resulting in a worse performance there.

%
%
\begin{table}[H]
\begin{center}
%
\tabcolsep=0.25cm    
\renewcommand*{\arraystretch}{1.15}              
\begin{tabular}{c|cccc}
                     & Theorem~\ref{CauchyTheorem} & Theorem~\ref{MatrixCauchyTheorem} & Theorem~\ref{TGC} (a) & Theorem~\ref{TGC} (b)      \\         
                     &    (Cauchy)                 &       (Matrix Cauchy)             &                       &                            \\   \hline
                     &                             &                                   &                       &                            \\
\text{Set 1}         &   1.26 / 8                  &      1.11 / 401                   &   1.11 / 529          &   1.13 / 62                \\   
                     &                             &                                   &                       &                            \\
\text{Set 2}         &   1.23 / 11                 &      1.07 / 200                   &   1.06 / 757          &   1.08 / 32                \\   
                     &                             &                                   &                       &                            \\
\text{Set 3}         &   1.58 / 0                  &      1.10 / 974                   &   1.15 / 17           &   1.11 / 9                 \\   
                     &                             &                                   &                       &                            \\
\text{Set 4}         &   1.48 / 0                  &      1.06 / 991                   &   1.10 / 3            &   1.07 / 6                 \\   
\end{tabular}
\caption{Cauchy-like results - ratio upper bound to modulus and number of times the bound outperformed the others.}  
\label{table1}        
\end{center}
\end{table}

\begin{table}[H]
\begin{center}
%
\tabcolsep=0.25cm    
\renewcommand*{\arraystretch}{1.15}              
\begin{tabular}{c|ccccc}
               & Theorem~\ref{PelletTheorem} & Theorem~\ref{MatrixPelletTheorem} & Theorem~\ref{TGP} (a) & Theorem~\ref{TGP} (b) & Theorem~\ref{TMGP}  \\        
               &   (Pellet)                  &   (Matrix Pellet)                 &                       &                       &                     \\  \hline
               &                             &                                   &                       &                       &                     \\
\text{Set 1}   &     1 / 0                   &             0                     &   0 / 0 / 203         &    0 / 0 / 145        &   0 / 0 / 214       \\   
               &                             &                                   &                       &                       &                     \\
\text{Set 2}   &   119 / 0                   &             0                     &   0 / 0 / 666         &    0 / 0 / 585        &   0 / 0 / 653       \\   
               &                             &                                   &                       &                       &                     \\
\text{Set 3}   &     0 / 0                   &            38                     &   4 / 0 / 345         &   17 / 0 / 243        &  38 / 0 / 0         \\   
               &                             &                                   &                       &                       &                     \\
\text{Set 4}   &     0 / 0                   &           976                     &  501 / 0  / 498       &    908 / 0 / 90       &  976/ 0 / 0         \\   
\end{tabular}
\caption{Pellet-like results - number of times inclusion regions for one and two zeros can be found.}  
\label{table2}        
\end{center}
\end{table}

\end{document}